\documentclass{article} 
\usepackage{amsmath}
\usepackage{amsthm}
\usepackage{amssymb}
\usepackage{mathrsfs}
\usepackage{paralist}
\usepackage{subcaption}
\usepackage[dvipsnames]{xcolor}
\usepackage[misc]{ifsym}
\usepackage{epsfig} 
\usepackage{epstopdf} 
\usepackage[colorlinks=true]{hyperref}
\hypersetup{urlcolor=blue, citecolor=red}
\allowdisplaybreaks

\usepackage[colorinlistoftodos]{todonotes}
\usepackage{mdframed}
\usepackage{xcolor}
\usepackage{bm}
\usepackage{changepage}
\usepackage{pdflscape} 
\usepackage{mathrsfs}
\usepackage[dvipsnames]{xcolor}
\usepackage{paralist}
\usepackage{booktabs}
\usepackage{bibentry}
\usepackage{float}
\usepackage{makecell}
\usepackage{algorithm}
\usepackage{algorithmic}
\usepackage{booktabs}
\usepackage{a4wide}
\usepackage{comment}
\usepackage[table]{xcolor}       
\usepackage{tablefootnote} 
\usepackage{stmaryrd}
\usepackage{makecell}
\usepackage{mathtools}
\usepackage{enumitem}
\usepackage{multirow}
\usepackage{pgfplots}
\pgfplotsset{compat=1.18}

\newtheorem{theorem}{Theorem}[section]

\theoremstyle{definition}

\newtheorem{remark}[theorem]{Remark}

\title{Physics-Informed Learning of Microvascular Flow Models using Graph Neural Networks}

\author{Paolo Botta$^{*,1}$, Piermario Vitullo$^{*,1}$, Thomas Ventimiglia$^{2}$, \\[1ex] Andreas Linninger$^{2,3}$ and Paolo Zunino$^{1,4}$}

\begin{document}

\maketitle

{\footnotesize
\centering
$^*$The first two authors equally contributed to the work \\
$^1$MOX, Department of Mathematics, Politecnico di Milano, Italy\\
$^2$Department of Bioengineering, University of Illinois at Chicago, IL, USA\\
$^3$Department of Neurosurgery, University of Illinois at Chicago, IL, USA\\
$^4$CMON Lab, Data Science Unit, Fondazione IRCCS Istituto Nazionale dei Tumori, Italy
\par} 



\begin{abstract}
The simulation of microcirculatory blood flow in realistic vascular architectures poses significant challenges due to the multiscale nature of the problem and the topological complexity of capillary networks. In this work, we propose a novel deep learning-based reduced-order modeling strategy, leveraging Graph Neural Networks (GNNs) trained on synthetic microvascular graphs to approximate hemodynamic quantities on anatomically realistic domains. Our method combines algorithms for synthetic vascular generation with a physics-informed training procedure that integrates graph topological information and local flow dynamics. To ensure the physical reliability of the learned surrogates, we incorporate a physics-informed loss functional derived from the governing equations, allowing enforcement of mass conservation and rheological constraints. The resulting GNN architecture demonstrates robust generalization capabilities across diverse network configurations. On the general-purpose test networks, the GNN reconstructs pressure, physical axial velocity recovered from the transformed model output, and hematocrit fields; on the cortical networks, the present quantitative validation is restricted to pressure. Although specific to the solver implementations and hardware used, the reported wall-clock comparisons show low-latency inference after training that may be useful in repeated-query workflows built on existing organ-scale graph-reduced solvers. This work establishes a new class of graph-based surrogate models for microvascular flow, grounded in physical laws and equipped with inductive biases that mirror mass conservation and rheological models, opening new directions for real-time inference in vascular modeling and biomedical applications.

\end{abstract}

\\




\def\eg{\textit{e.g.}\ }
\def\ie{\textit{i.e.}\ }
\def\argmin{\mathrm{argmin}}
\def\esssup{\mathrm{ess\,sup}}
\def\I{\boldsymbol{\mathrm{I}}}
\def\Tau{\mathcal{T}}
\def\para{\boldsymbol{\mu}}
\def\parmi{\boldsymbol{\mu}_{m}}
\def\parma{\boldsymbol{\mu}_{M}}
\def\micro{m}
\def\macro{M}
\def\parph{\boldsymbol{\mu}_{ph}}
\def\parg{\boldsymbol{\mu}_{g}}
\def\domega{\partial \Omega}
\def\uph{u_{\para,h}}
\def\uuphi{\mathtt{\bf u}_{\para^{(i)},h}}
\def\uuh{\mathtt{\bf u}_{h}}
\def\uuph{\mathtt{\bf u}_{\para,h}}
\def\uphrb{u_{\para,h}^{rb}}
\def\nuphrb{\widetilde{u}_{\para,h}^{rb}}
\def\nuuphrb{\widetilde{\bf u}_{\para,h}^{rb}}
\def\muphrb{\widehat{u}_{\para,h}^{rb}}
\def\muuphrb{\widehat{\mathbf{u}}_{\para,h}^{rb}}
\def\muuphrbj{\widehat{\mathbf{u}}_{\para,h}^{rb,(j)}}
\def\uuphrb{\mathtt{\bf u}_{\para,h}^{rb}}
\def\urb{\mathtt{u}_{\para,rb}}
\def\nurb{\widetilde{\mathtt{u}}_{\para,rb}}
\def\nuu{\widetilde{\mathtt{u}}_{\para}}
\def\nuurb{\widetilde{\mathtt{\bf u}}_{\para,rb}}
\def\nuuu{\widetilde{\mathtt{\bf u}}_{\para,c}}
\def\muurb{\widehat{\mathtt{\bf u}}_{\para,rb}}
\def\uurb{\mathtt{\bf u}_{\para,rb}}
\def\npmicro{n_{\micro}}
\def\npmacro{n_{\macro}}
\def\nrb{n_{rb}}
\def\nnrb{\widetilde{n}_{rb}}
\def\npsi{\widetilde{\psi}}
\def\nft{\widetilde{\alpha}}
\def\RR{\mathbb{R}}
\def\VV{\mathbb{V}}
\def\WW{\mathbb{W}}
\def\SS{\mathbb{S}}
\def\ann{\mathcal{N}}
\def\lnn{\mathcal{L}}
\def\minn{\mathcal{M}}
\def\hyper{\boldsymbol{\theta}}
\def\hypermi{\boldsymbol{\theta}_{\micro}}
\def\hyperd{\boldsymbol{\theta}_{d}}
\def\hyperi{\boldsymbol{\theta}_{\eta}}
\def\hyperma{\boldsymbol{\theta}_{\macro}}
\def\wwmic{\boldsymbol{\widetilde{\theta}}_{g}}
\def\uu{\mathbf{v}}
\def\inlets{\boldsymbol{\eta}}
\def\dist{\mathbf{d}}
\newcommand{\bct}{{C}_t}
\newcommand{\bp}{{p}}
\newcommand{\argminF}{\mathop{\mathrm{argmin}}\limits}

\def\weight{\boldsymbol w}
\def\feature{\widetilde{\boldsymbol{w}}}

\newcommand{\rev}[1]{{\color{teal} #1}}


\section{Introduction}
\label{sec:introduction}
The simulation of microcirculatory blood flow is an essential component in understanding tissue perfusion, vascular function, and the development of pathological conditions such as ischemia, inflammation, and cancer~\cite{Popel2005, Secomb2017a,Figueroa2017}, just to give a few examples. On the microscopic scale, the transport of blood through networks of arterioles, capillaries, and venules is governed by intricate hemodynamic phenomena and highly complex topological structures, which pose significant computational challenges.

High-fidelity models are able to accurately capture these dynamics. However, their applicability to large-scale, anatomically realistic microvascular networks is severely hampered by the associated computational burden and the difficulties in accommodating topological variability~\cite{linninger1,linninger2,linningerplos}. 
In the case of fully resolved three-dimensional models (including or not the description of red blood cells) even in relatively small tissue volumes, such as a cortical sample of size $1 \, mm^3$~\cite{linninger2,GL2,Schmid2017,Schmid2019}, the resolution required to represent individual vessels can easily lead to systems not yet computable using a direct numerical simulation approach. 
To overcome this limitation, dimensional model reduction has been adopted with the purpose to reduce the mathematical and computational complexity of blood flow to one-dimensional models defined along the vessel centerline, such as in the \textit{geometric multiscale methods} \cite{QUARTERONI2016193}. This approach can be modulated to a high level of physiological and anatomical detail~\cite{alastruey2011pulse,muller2016high,Possenti2019a,ventimiglia2023meshfree,doi:10.1177/0271678X231214840}, despite the intrinsic multiscale nature of microvascular flow~\cite{karst2015oscillations,Penta2015,Peyrounette2018,Tong2021}. State of the art simulations reach today the organ scale level~\cite{doi:10.1177/0271678X231214840,10.1371/journal.pcbi.1013459}. For an accurate description of the pulsatile regime such simulations still require high performance computing resources, while the computational cost is substantially reduced in the case of the steady regime typical of the microcirculation. 
Despite the great success of one-dimensional model for the simulation of blood flow, the research on surrogate models for blood flow in complex networks is still active~\cite{Pfaller2024449,VITULLO2024104068}, with the aim to preserve predictive accuracy while achieving substantial computational speed-ups. These advancement will extend the ability to simulate pulsatile blood flow beyond the current limits and also enable the use of the simpler microcirculation models in the context of computational workflows involving many-query simulations such as inverse problem formulations, optimization, and uncertainty quantification
Vascular networks admit a natural representation as graphs, with vessels corresponding to edges and junctions to nodes. GNNs are designed to operate directly on such graph-structured inputs, learning approximations of the underlying physics from data~\cite{Bronstein2017,hamilton2017inductive,battaglia2018relational}.
Moreover, GNNs can be trained in a physics-informed fashion, incorporating domain-specific knowledge (such as mass conservation and nonlinear rheology) directly into the loss functional. This strategy mitigates the typical limitations of purely data-driven models and enhances the ability to generalize to unseen network geometries. 
After the offline training stage, GNN inference consists of a forward evaluation and can provide low-latency repeated predictions. The magnitude of any wall-clock advantage depends on the compared solver implementations and hardware. For the steady models studied here, a single full-order solve is already inexpensive. The principal computational motivation is therefore repeated evaluation in many-query settings, rather than making a previously intractable solve possible.

In this work, we develop a physics-informed computational learning problem specifically designed to address flow problems on graphs using GNN architectures. We apply this framework to define surrogate models for microcirculatory blood flow in the steady regime, although ongoing work confirms that the approach can be successfully extended to the pulsatile regime \cite{behrens2025}. 
We show that a model trained exclusively on synthetic capillary networks transfers its pressure predictions to large-scale, anatomically realistic vascular reconstructions.
Remarkably, we show that a model trained exclusively on synthetic capillary networks with generic topological and geometric features is able to generalize robustly to large-scale, anatomically realistic vascular reconstructions. 

In particular, when tested on mouse cerebral cortex networks that exhibit a substantially higher degree of complexity than training data \cite{SCHREINER199927,linninger1,linninger2,linningerplos}, the GNN maintains high predictive precision on the blood pressure, reproducing full-order simulation results with excellent agreement. Starting from encouraging insight and results, accurate prediction of cortical velocity, hematocrit, and RBC-flux conservation remains future work.

Compared to previous studies on GNNs for blood flow \cite{Pegolotti2024,Pfaller2024449,Iacovelli2024}, our approach introduces several differences. Existing work, such as the reduced-order models of \cite{Pegolotti2024} and their recent extensions that integrate recurrent architectures for temporal learning \cite{Iacovelli2024}, focuses primarily on large-scale cardiovascular flow within arterial trees and exploits data-driven surrogates of one-dimensional flow dynamics. In contrast, we address the steady microcirculatory regime on the scale of large vascular networks, where flow is governed by mass conservation, diffusion mechanisms, and nonlinear rheology. For extensions to the pulsatile regime, we refer to Section \ref{sec:conclusions}.

The remainder of this manuscript is organized as follows. Section \ref{sec:method} introduces the general mathematical framework for physics-informed learning on graphs, formalizing the parameter-to-solution map and its approximation via GNNs. Section \ref{sec:models} presents the physical models for microvascular blood flow, encompassing both linear and nonlinear rheology, while Section \ref{sec:networks} details the algorithms used to generate synthetic and anatomically realistic vascular networks. Section \ref{sec:gnn_models} describes the architecture, loss formulations, and training strategies of the proposed GNN surrogates, highlighting the integration of physical constraints within the learning process. Finally, Section \ref{sec:results} reports the numerical experiments that evaluate the predictive accuracy and generalization of the models in synthetic and image-based microvascular networks, followed by the concluding remarks.

\section{Physics informed learning using graph neural networks}
\label{sec:method}

In this section, we introduce the methodological framework underlying the proposed approach. We formulate parameterized problems on graphs, providing a general setting for network-based physical systems that couple discrete topological structures with continuous quantities defined on nodes and edges. The associated parameter-to-solution map thus exhibits a mixed discrete–continuous nature. Graph Neural Networks (GNNs) are naturally suited to fit this structure, leveraging graph connectivity to propagate information and learn nonlinear relations between nodal and edge variables. This formulation yields a unified paradigm for the operator approximation on graphs, ideal for the representation of microvascular blood flow in complex configurations.

\subsection{A learning problem on graphs}
\label{sec:abstract_metric_graph_problem}


Let $G := (V, E)$ be a finite, \emph{directed graph} encoding the topological connectivity of the network, where $V := \{ v_j \}_{j=1}^{n}$ denotes the set of $n$ vertices in $\mathbb{R}^3$ and $E := \{ e_i \}_{i=1}^{m}$ is the set of $m$ directed edges that connect pairs of vertices. In particular, each edge $e = (v_{s}, v_{t}) \subset \mathbb{R}^3$ sets the orientation from the source node $v_{s}$ to the target node $v_{t}$.
All relevant information on the connectivity of a directed graph is encoded in the connectivity (or incidence) matrix $C \in \mathbb{Z}^{m\times n}$ with entries:
\[
    C_{ij} := \{
        1, \ \text{if } e_i \text{ enters } v_j; \quad
       -1, \ \text{if } e_i \text{ leaves } v_j; \quad
        0, \ \text{otherwise}\}.
\]

We then introduce a vector of weights defined on the vertices of $G$, $W^V :=\{\weight^V_j \}_{j=1}^{n}$, and on the edges, $W^E :=\{\weight^E_i \}_{i=1}^{m}$, where $\weight^V_j \in \mathbb{R}^{p_V}$ $\forall$ $j=1,...,n$ and $\weight^E_i \in \mathbb{R}^{p_E}$ $\forall$ $i=1,...,m$. This leads to the definition of the vector of weights $W := 
\big[W^V, W^E\big]^{\!\top} \in \mathbb{R}^{n\,p_V + m\,p_E}$.
A \emph{weighted directed graph} is the triple $\mathcal{G} := (V, E, W)$.
Each graph (directed or weighted) $\mathcal{G}$ can be embedded in $\mathbb{R}^3$, proceeding as follows:
\begin{description}
    \item[\textit{(i)}] we associate with each node  in $V$ a point  $\mathbf{x}_i \in \mathbb{R}^d$ in a bijective way, with $d$ denoting the number of spatial dimensions ($d = 2,3$). Consequently, we introduce the matrix that encodes the coordinates of each vertex, $X \in \mathbb{R}^{n \times 3}$ such that $X(i,:)= \mathbf{x}_i^\top$;
    \item[\textit{(ii)}] we connect two points $\mathbf{x}_v, \mathbf{x}_w$ by a three-dimensional simple arc $\mathbf{s}_{vw}$ if and only if $v, w$ are adjacent, in such a way that different arcs do not share any internal points;
    \item[\textit{(iii)}] we define a length function $\rho : E \to \mathbb{R}_+$ assigning a metric on each edge $e = (v_{s}, v_{t})$, such that its length is $L = \rho(e)$.
\end{description}
We know that for any weighted directed graph $\mathcal{G}$, the resulting embedded object in $\mathbb{R}^3$ is a manifold. In particular, assuming for all edges $e$ the convention of identifying the source node $v_s$ with 0 and the target node $v_t$ with $\rho(e)$, it is possible to define a metric over $\mathcal{G}$. 

This construction is the root of the notion of a \emph{metric graph}. Let $\mathcal{G}=(V, E, W)$ be a weighted directed graph and let us define
$\mathcal{E} := \prod_{e \in E} \{ e \} \times (0, \rho(e))$.
Then the triple $\mathfrak{G} := (V, \mathcal{E}, W)$
is called the (weighted) \emph{metric graph} over $\mathcal{G}$.
A point of $\mathfrak{G}$ is an element of either $V$ or $\mathcal{E}$. The elements of $\mathcal{E}$ are called metric edges.

We now consider steady-state problems governed by differential operators acting on $\mathfrak{G}$.
Under these assumptions, the most general problem we can consider is of the form
\begin{equation}\label{eq:metric_graph_problem}
    \begin{dcases}
        \mathfrak{L}^W u^W= f^W 
        & \text{in } \mathfrak{G}, \\
        \mathfrak{B}^W u^W = g^W
        & \text{on } \partial \mathfrak{G},
    \end{dcases}
\end{equation}
where $\mathfrak{L}^W$ denotes an operator on the metric graph, and $\mathfrak{B}^W$ represents boundary or vertex conditions enforcing continuity and conservation properties across the nodes of $\mathfrak{G}$.
A prototypical example of $\mathfrak{L}^W$ is a possibly semilinear elliptic partial differential equation in $\mathfrak{G}$, where the weights $W$ represent the parameters of the operator $\mathfrak{L}$ on the edges $[0,L_i]$.

However, as the physical models used for microvascular blood flow are rooted in Poiseuille's law, the problems addressed here will reduce to nonlinear algebraic equations on the weighted directed graph $\mathcal{G}=(V,E,W)$. 
Let us consider for example the following linear relations:
\begin{equation}\label{eq.general_poiseuille}
    \left\{
    \begin{aligned}
    u^V_{s(i)} - u^V_{t(i)} &= R_i^E u^E_{i}, & \quad i &= 1,\ldots,m,\\[2mm]
    \sum_{i\in E(v_j)} C_{ij} u^E_{i} &= 0, 
    & \quad j &= 1,\ldots,n.
    \end{aligned}
    \right.
\end{equation}
where $C$ is the incidence matrix, and $R^E_i$ are the edge resistances.
Setting $\mathbf{u}_{j}^{V}\in\mathbb{R}^{p_V'}$ $\forall$ $j=1,...,n$ and $\mathbf{u}_{i}^{E}\in\mathbb{R}^{p_E'}$  $\forall$ $i=1,...,m$, this problem can be formulated in terms of unknown vector values on the vertices and edges of the graph, defined as follows:
\[
\mathbf{U}^V := \{ \mathbf{u}_{j}^{V} \}_{j=1}^{n},
\quad
\mathbf{U}^E := \{ \mathbf{u}_{i}^{E} \}_{i=1}^{m},
\quad\mathbf{U}^W :=
\begin{bmatrix}
\mathbf{U}^V \\[4pt]
\mathbf{U}^E
\end{bmatrix}.
\]
As a result, problem \eqref{eq.general_poiseuille} becomes the following:
\begin{equation}\label{eq:discrete_system}
    L^W \, \mathbf{U}^W = \mathbf{F}^W,
    \quad 
    L^W :=
    \begin{bmatrix}
        C & - R^E \\[10pt]
        0 & C^{\!\top}
    \end{bmatrix},
    \quad
    \mathbf{F}^W:=
    \begin{bmatrix}
        \mathbf{F}^V \\[10pt]
        \mathbf{F}^E 
    \end{bmatrix},
\end{equation}
 where $R^E_{i}=r(\weight_i)$, being $r: \mathbb{R}^{p_V+p_E} \rightarrow \mathbb{R}$ a suitable constitutive function. The right-hand side $\mathbf{F}^W$ encodes the boundary conditions that will be addressed later.
The system can be reduced to the vertex unknowns, with the corresponding matrix $C^{\!\top} (R^E)^{-1} C$ defining a graph Laplacian, which is a possible finite-dimensional realization of the continuous operator $\mathfrak{L}^W$.
The solution depends on both the weights $W$ (representing the physical parameters) and the (simple or directed) graph $G$, encoded by the matrix of vertices $X \in \mathbb{R}^{n \times 3}$ and the discrete incidence matrix $C \in \mathbb{Z}^{m\times n}$. We therefore define the \emph{parameter-to-solution map}
\[
\mathcal{S} : \mathbb{R}^{n\times d} \times \mathbb{Z}^{m\times n} \times \mathbb{R}^{n\,p_V + m\,p_E}
\longrightarrow 
\mathbb{R}^{n\,p_V' + m\,p_E'},
\quad 
\mathcal{S}(G,W) = \mathbf{U}^W.
\]
Given $\widehat{\mathcal{H}}$ a set of candidate functions for this map, $\widehat{\mathcal{S}} \in \widehat{\mathcal{H}}$, given a collection of labeled data $\mathcal{S}(G_k,W_k)=\mathbf{U}_k^W$, for $k=1,\ldots,N$, given a loss function measuring the discrepancy between data and predictions, $l(\widehat{\mathbf{U}}^W,\mathbf{U}^W)$,  and a feature mapping of the weight vector $\widetilde{W}:=\Phi(W)$, we aim to solve the following supervised learning problem, that consists of the minimization of the empirical error $\mathcal{L}$ (from now on called \textit{empirical total loss} or just \textit{loss}, with little abuse of notation) over the candidate functions:
\begin{align*}
    &\mathcal{L}\big((\widehat{\mathcal{S}}(G,\widetilde{W}),\mathbf{U}^W);(G_k,W_k),\mathbf{U}_k^W\big)=\frac{1}{N}\sum_{k=1}^Nl(\widehat{\mathbf{U}}_k^W,\mathbf{U}_k^W),
\\
    &\widehat{\mathcal{S}}^* = \argminF_{\widehat{\mathcal{S}} \in \widehat{\mathcal{H}}} \, \mathcal{L}\big((\widehat{\mathcal{S}}(G,\widetilde{W}),\mathbf{U}^W);(G_k,W_k),\mathbf{U}_k^W\big).
\end{align*}
To solve this learning problem, we adopt a Graph Neural Network architecture, which is ideally suited to approximate a parameter-to-solution map featuring a mixed mathematical structure combining discrete and continuous information.
In fact, GNNs~\cite{Scarselli2009,Bronstein2017,battaglia2018relational} process node- and edge-based continuous features through standard feed-forward neural networks, while a message passing mechanism exploits the discrete incidence structure of the underlying graph to propagate information across connected entities. 

\subsection{Graph Neural Networks architecture}
\label{sec:GNNA}

Building upon the formalism introduced in Section~\ref{sec:abstract_metric_graph_problem}, we recall that the physical domain is a weighted directed graph $\mathcal{G} = (V, E, W)$. For computational purposes, we consider its discrete counterpart  $G = (V, E)$, which retains the same topological connectivity, while the physical and geometric quantities are encoded within the node and edge weights $W = [W^V, W^E]^{\!\top}$.

Within this setting, we adopt a \emph{supervised learning} approach to approximate the parameter-to-solution map introduced in Eq.~\eqref{eq:discrete_system}, acting on weighted directed graphs. Formally, the goal is to learn an approximation 
$\widehat{\mathcal{S}}$ of the exact solution operator $\mathcal{S}$, namely:
\[
\widehat{\mathbf{U}}^{W} = 
\widehat{\mathcal{S}}(G, \widetilde{W}) 
\approx 
\mathcal{S}(G, W) = \mathbf{U}^{W}.
\]

\begin{figure}[t]
    \centering
    \includegraphics[width=1\textwidth]{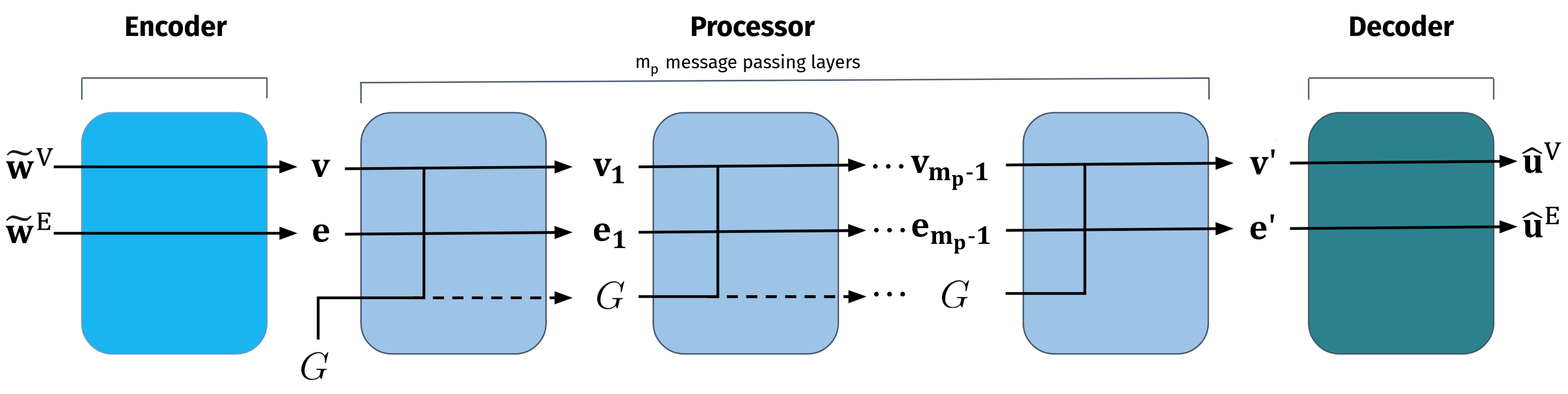}
    \caption[GNN complete architecture]{GNN complete architecture, consisting of an encoder, a processor, and a decoder.}
    \label{fig:gnn_architecture}
\end{figure}

The adopted GNN architecture follows the standard \textit{Encoder--Processor--Decoder} paradigm, schematically illustrated in Figure~\ref{fig:gnn_architecture}. 
The \textit{encoder} (denoted with $\mathcal{E}_v(\feature^V),\,\mathcal{E}_e(\feature^E)$) for the vertices and edges, respectively) maps the physical and geometric input features defined on the nodes and edges into a latent space of fixed dimension. 
The \emph{processor}, composed of multiple message passing layers (named $\mathcal{P}_v(G, \mathcal{E}_v(\feature^V), \mathcal{E}_e(\feature^E))$ and $\mathcal{P}_e(G, \mathcal{E}_v(\feature^V), \mathcal{E}_e(\feature^E))$), iteratively propagates and updates information throughout the graph according to its incidence structure. We note that all the message passing layers share a fixed hidden dimension $l$, updating node features $\mathbf{v} : V \to \mathbb{R}^l$ and edge features $\mathbf{e} : E \to \mathbb{R}^l$ without modifying their dimension. Finally, \emph{decoders} $\mathcal{D}_v(\mathcal{P}_v(G, \mathcal{E}_v(\feature^V), \mathcal{E}_e(\feature^E)))$ and $\mathcal{D}_e(\mathcal{P}_e(G, \mathcal{E}_v(\feature^V), \mathcal{E}_e(\feature^E)))$ transform the latent representations of the nodes and edges into the predicted physical quantities. 
More precisely, the composition of these maps can be represented as follows:
\begin{equation*}
 \widehat{\mathbf{u}}^V = \mathcal{D}_v\Big(\mathcal{P}_v\big(G, \mathcal{E}_v(\feature^V), \mathcal{E}_e(\feature^E)\big)\Big), \qquad
 \widehat{\mathbf{u}}^E = \mathcal{D}_e\Big(\mathcal{P}_e\big(G, \mathcal{E}_v(\feature^V), \mathcal{E}_e(\feature^E)\big)\Big).
\end{equation*}
The resulting outputs are assembled into the vector-valued approximation of the solution, 
which provides the learned estimate of the parameter-to-solution map:
\[
\widehat{\mathcal{S}}(G, \widetilde{W}) = \widehat{\mathbf{U}}^W :=
\begin{bmatrix}
\widehat{\mathbf{U}}^V \\[4pt]
\widehat{\mathbf{U}}^E
\end{bmatrix}, \quad
\widehat{\mathbf{U}}^V := \{ \widehat{\mathbf{u}}_{j}^{V} \}_{j=1}^{n},
\quad
\widehat{\mathbf{U}}^E := \{ \widehat{\mathbf{u}}_{i}^{E} \}_{i=1}^{m}
.
\] 
This architecture, illustrated in Figure~\ref{fig:gnn_architecture}, enables simultaneous treatment of nodal and edge quantities in the processor and in the decoder, which is essential for predicting coupled variables such as pressure and flow.

\subsubsection{Feedforward Neural Networks}
Feedforward Neural Networks (FFNNs) constitute the fundamental building blocks of most modern deep learning architectures. They are composed of a sequence of layers that apply affine transformations followed by nonlinear activation functions, enabling the network to approximate complex nonlinear mappings between input and output spaces.  
Given an input vector $\mathbf{x} \in \mathbb{R}^{d_{\text{in}}}$, a generic FFNN with $K$ layers can be expressed as the composition:
\[
\Psi(\mathbf{x};\theta) = \mathfrak{f}_K \circ \mathfrak{f}_{K-1} \circ \dots \circ \mathfrak{f}_1(\mathbf{x}),
\quad
\mathfrak{f}_\ell(\mathbf{z}) = \sigma_\ell(\mathcal{W}_\ell \mathbf{z} + \mathbf{b}_\ell),
\]
where $\mathcal{W}_\ell$ and $\mathbf{b}_\ell$ denote the weight matrix and bias vector of the $\ell$-th layer, respectively, $\sigma_\ell$ is a nonlinear activation function (e.g., ReLU, $\tanh$), and $\theta = \{\mathcal{W}_\ell, \mathbf{b}_\ell\}_{\ell=1}^{K}\in\Theta$ represents the set of all trainable parameters of the network.  

In the present framework, FFNNs are employed throughout the architecture wherever a nonlinear map $\Psi$ appears.   
For simplicity of notation, we omit the explicit dependence on the parameters $\theta$ in the subsequent sections, writing simply $\Psi(\cdot)$ instead of $\Psi(\cdot;\theta)$, with the understanding that all feedforward mappings implicitly depend on their respective trainable weights and biases.

\subsubsection{The Encoder and Decoder modules}
Each node $v \in V$ and edge $e \in E$ is associated with feature vectors that represent local physical or geometric information. \textit{Node input features} $\widetilde{W}^V := \{\feature^V_j\}_{j=1}^{n}$, with $\feature^V_j  \in \mathbb{R}^{\widetilde{p}_V}$, encode local physical quantities defined at each vertex (e.g., boundary conditions, nodal parameters). \textit{Edge input features} $\widetilde{W}^E := \{\feature^E_i\}_{i=1}^{m}$, with $\feature^E_i  \in \mathbb{R}^{\widetilde{p}_E}$, encode the geometric or material properties associated with each edge (e.g., length, diameter, resistance). The specific form of $\feature^V$ and $\feature^E$ is problem-dependent. The collection of these features is $\widetilde{W}:=[\widetilde{W}^V, \widetilde{W}^E]^{\!\top} := \Phi (W)$ as a feature mapping of the weight vector. We stress that the mapping $\Phi$ is a hyperparameter specified \emph{a priori} based on the problem structure. Its choice determines which weights and aspects of the physical system are exposed to the learning architecture. Since the dimensions of the input features generally differ from $l$, we employ feedforward neural networks to map them to the latent space:
\begin{equation*}
    \mathcal{E}_v(\feature^V) = \Psi_{\mathcal{E}}^v(\feature^V)
    \qquad
    \mathcal{E}_e(\feature^E) = \Psi_{\mathcal{E}}^e(\feature^E)
\end{equation*}
The GNN output features corresponding to node- and edge-level quantities are the \textit{node output features} $\mathbf{U}^V := \{\mathbf{u}^V_j\}_{j=1}^{n}$, with $\mathbf{u}^V_j \in \mathbb{R}^{p_V'}$, representing physical states at the nodes (e.g., pressure). \textit{Edge output features} $\mathbf{U}^E := \{\mathbf{u}^E_i\}_{i=1}^{m}$, with $\mathbf{u}^E_i \in \mathbb{R}^{p_E'}$, correspond to the quantities defined along the edges (for example, flow rate). 

The processor module outputs hidden node and edge features, $\mathbf{v'} = \mathcal{P}_v(G, \mathbf{v,e})$ and $\mathbf{e'} = \mathcal{P}_e(G,\mathbf{v,e})$. 
Since these do not match the output dimensionalities, two separate feedforward neural networks are employed to recover the physical output features:
\begin{equation*}  
    \mathcal{D}_v(\mathbf{v'}) = \Psi_{\mathcal{D}}^v(\mathbf{v'}), \qquad
    \mathcal{D}_e(\mathbf{e'}) = \Psi_{\mathcal{D}}^e(\mathbf{e'}).
\end{equation*}

While the terms \textit{encoder} and \textit{decoder} are traditionally associated with dimensionality reduction and reconstruction, in our context they are used more broadly. The encoder extracts latent features from the input, and the decoder maps latent representations to physically meaningful outputs. Depending on the task, these modules may either increase or decrease the dimension of the features.



\subsubsection{The Processor module}
The processor module is composed of $m_p$ message passing steps that iteratively update both the node and the edge features. Let $H_j$ and $F_j$ denote the edge and node update functions. Starting from encoded features, the recursive update now reads:
\begin{equation*}
    \begin{cases}
    \mathbf{v}_0 = \mathcal{E}_v(\feature^V), \quad \mathbf{e}_0 = \mathcal{E}_e(\feature^E), & \\
    \mathbf{e_j} = H_j(G, \mathbf{v}_{j-1},\mathbf{e}_{j-1}) & j=1,\dots,m_p, \\
    \mathbf{v_j} = F_j(G, \mathbf{v}_{j-1},\mathbf{e}_j) & j=1,\dots,m_p,
    \end{cases}
\end{equation*}
so that the processor outputs
\begin{equation*}
    \mathcal{P}_v(G, \mathbf{v,e}) = F_{m_p}(G, \mathbf{v}_{m_p-1},\mathbf{e}_{m_p}), \quad
    \mathcal{P}_e(G, \mathbf{v,e}) = H_{m_p}(G, \mathbf{v}_{m_p-1},\mathbf{e}_{m_p-1}).
\end{equation*}
This formulation makes explicit that the node and edge states evolve in tandem at each message passing step, by iteratively applying the message passing layers $m_p$. 
The parameter $m_p$ is a hyperparameter of the architecture: a processor with $m_p$ layers enables communication between nodes that are $m_p$ edges apart, thereby extending the receptive field from local to increasingly nonlocal interactions.

The key component of this algorithm is the \textit{message passing layer}, which propagates information along graph edges~\cite{Scarselli2009,Bronstein2017,battaglia2018relational}. In its standard form, the update acts only on node features, while edge features remain unchanged. Specifically, given hidden node features $\mathbf{v}$ and hidden edge features $\mathbf{e}$, the message passing layer $F=F(G, \mathbf{v},\mathbf{e})$ takes as input $\mathbf{v}$, $\mathbf{e}$ and $G$ and outputs a new collection of vertex features $\mathbf{v}' : V \to \mathbb{R}^l$.

In order to define nodal update function $F$, let the mappings $\mathbf{v} : V \to \mathbb{R}^l$ and $\mathbf{e} : E \to \mathbb{R}^l$ denote hidden node and edge features of dimension $l$. 
For a directed edge $(v_s, v_t) \in E$, we set
\[
\mathbf{v}_{\text{in}}(v_s, v_t) = \mathbf{v}(v_s), \quad
\mathbf{v}_{\text{out}}(v_s, v_t) = \mathbf{v}(v_t).
\]
In particular, $\mathbf{v}_{\text{in}}$ computes the vertex features of the corresponding source node, while $\mathbf{v}_{\text{out}}$ computes the vertex features of the corresponding destination node.
The aggregation function for a node $v \in V$ is then defined as
\[
\overline{\mathbf{e}}(v) = \sum \mathbf{e}(v_s, v),
\ \forall (v_s, v) \in E, \ v_s \in V,
\]
which collects the contributions of all incoming edges.
Given two functions with the same domain $\mathbf{f}:X \to \mathbb{R}^a$ and $\mathbf{g}:X \to \mathbb{R}^b$, the concatenation operator $\oplus$ is defined as
\[
(\mathbf{f} \oplus \mathbf{g})(x) = [f_1(x),...,f_a(x),g_1(x),...,g_b(x)]
\]
where $x \in X$ is a generic input.
The standard message passing update reads
\[
F(G, \mathbf{v}, \mathbf{e}) 
= \Psi_u \left( \mathbf{v} \oplus \Psi_M \left(\overline{ \mathbf{e} \oplus \mathbf{v}_{\text{in}} \oplus \mathbf{v}_{\text{out}} }\right) \right),
\]
which concatenates and aggregates contributions from incoming edges, and $\Psi_u$ and $\Psi_M$ are feedforward neural networks, such that
$\Psi_u:\mathbb{R}^{2l} \to \mathbb{R}^l$ and $\Psi_M:\mathbb{R}^{3l} \to \mathbb{R}^l$.
In GNN architectures, $F$ is considered a transformer that modifies the characteristics of the vertices connected to the graph nodes. A message passing layer, within these architectures, refers to a specific type of graph-forwarding routine, which exploits the local structure of the graph $G$. This routine allows communication exclusively between neighboring nodes. 

\begin{figure}[t]
    \centering
    \includegraphics[width=\textwidth]{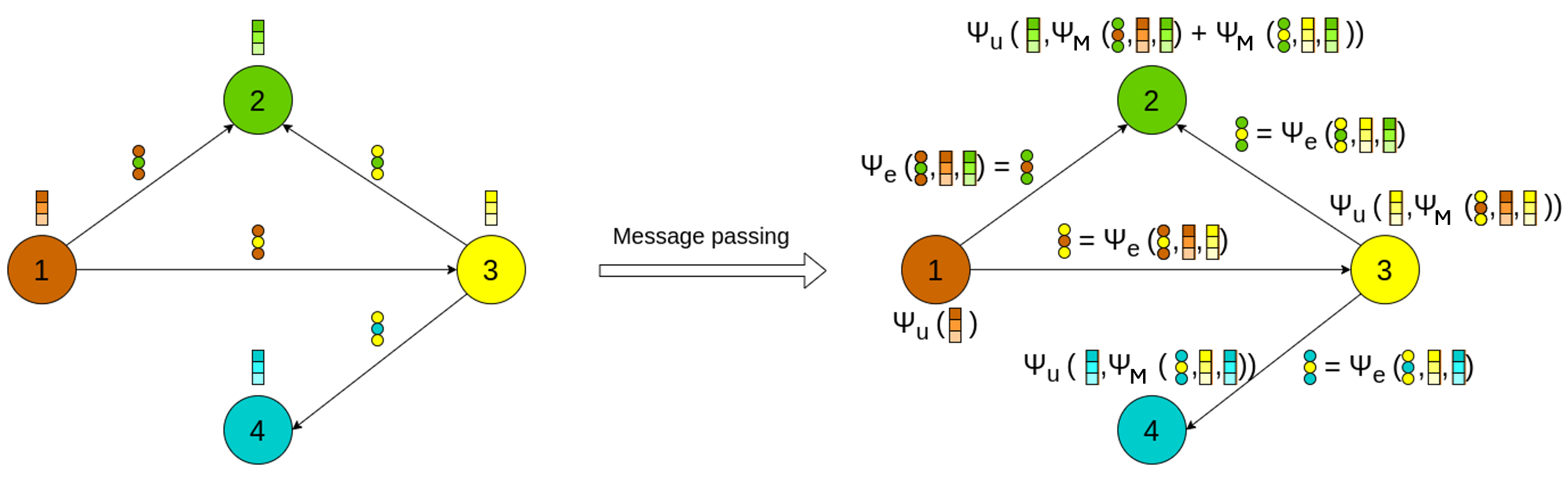}
    \caption{Message propagation and aggregation with edge features update.}
    \label{fig:mps2}
\end{figure}

In what follows, we will adopt an extended message passing mechanism that updates both node and edge states through local interactions. Specifically, the update rules are:
\[
H(G, \mathbf{v}, \mathbf{e}) = \Psi_e \left( \mathbf{e} \oplus \mathbf{v}_{\text{in}} \oplus \mathbf{v}_{\text{out}} \right), 
\quad
F(G, \mathbf{v}, \mathbf{e}) = \Psi_u \left( \mathbf{v} \oplus \Psi_M \left(\overline{ H \oplus \mathbf{v}_{\text{in}} \oplus \mathbf{v}_{\text{out}}}\right) \right),
\]
where $\Psi_e:\mathbb{R}^{3l} \to \mathbb{R}^l$ is an additional feedforward neural network.
Note that the nonlinearity of the message passing update arises from the activation functions within the neural networks $\Psi_e, \Psi_M, \Psi_u$. 
This formulation makes explicit that both edge and node features evolve at each iteration, a critical property when modeling vascular networks where edges carry hemodynamic variables (e.g., flow, diameter) and nodes represent pressures or bifurcation states.
For a graphical representation of this new strategy, see Figure~\ref{fig:mps2}. 

To improve gradient flow and training stability, we include skip (residual) connections between consecutive message passing layers. This strategy, originally introduced in ResNet architectures~\cite{he2016deep}, has proven effective in GNNs as well~\cite{xu2021optimization}. 
Concretely, the processor module with $k>0$ skip connections becomes:
\begin{equation*}
    \begin{cases}
            \mathbf{v}_0 = \mathcal{E}_v(\feature^V), \quad \mathbf{e}_0 = \mathcal{E}_e(\feature^E), \\
    \mathbf{e_j} = H_j(G, \mathbf{v}_{j-1},\mathbf{e}_{j-1}) + \mathbf{e}_{j-k} \, \mathbb{I}(j \geq k) & j=1,\dots,m_p, \\
    \mathbf{v_j} = F_j(G, \mathbf{v}_{j-1},\mathbf{e}_j) + \mathbf{v}_{j-k} \, \mathbb{I}(j \geq k) & j=1,\dots,m_p,
    \end{cases}
\end{equation*}
Skip connections are particularly beneficial when stacking a large number of message passing layers ($m_p \gg 1$), as they allow information to propagate across the network without degradation, enabling more expressive but still trainable architectures.


Finally, although the vascular network naturally defines a directed graph, reflecting the physical orientation of blood flow, we model the input graph provided to the GNN as undirected. This choice allows the message-passing scheme to propagate information bidirectionally along each vessel segment, enabling both source-to-target and target-to-source interactions at every iteration. Empirically, we observe that this symmetric exchange of information significantly improves convergence of the training algorithm and predictive accuracy.

\subsection{Loss functions and training}
\label{sec:loss_function}


To train the GNN surrogate models introduced in Section~\ref{sec:GNNA}, we adopt composite loss functions that incorporate both data fidelity, as in the supervised learning setting, and problem-specific physical information within a physics-informed framework. As a consequence, the total loss function used to train the GNN is defined as:
\begin{equation}
    \mathcal{L}(\widehat{\mathcal{S}},\mathbf{U}^W; \theta) := \delta \, \mathcal{L}_{\text{data}}(\theta) + (1 - \delta) \, \mathcal{L}_{\text{physics}}(\theta),
    \label{eq:total_loss}
\end{equation}
where $\mathcal{L}_{\text{data}}$ ensures data fidelity and $\mathcal{L}_{\text{physics}}$ is the physics-informed component.
The symbol $\theta$ denotes the trainable parameters of the network, and $\delta \in [0,1]$ is a tunable scalar that weights the relative importance of the data-driven and physics-informed components.

Unless otherwise specified, from now on we denote the Euclidean norm of the vector $\mathbf{v}$ by $\|\mathbf{v}\|$. This is also called the $L^2$ norm of the vector.
The data-fidelity loss quantifies the mismatch between the GNN predictions and the full-order model (FOM) ground truth on a dataset of $N$ training graphs:
\begin{equation}
    \mathcal{L}_{\text{data}}(\theta) = \frac{1}{N} \sum_{k=1}^N \left( 
    \frac{1}{n_k} \sum_{j = 1}^{p'_V} \alpha_j \|\mathbf{u}^V_{j,k} - \widehat{\mathbf{u}}^V_{j,k}(\theta)\|^2 + 
    \frac{1}{m_k} \sum_{i = 1}^{p'_E} \beta_i \|\mathbf{u}^E_{i,k} - \widehat{\mathbf{u}}^E_{i,k}(\theta)\|^2 
    \right),
    \label{eq:data_loss}
\end{equation}
where $\{\alpha_j\}_{j=1}^{p'_V}, \{\beta_i\}_{i=1}^{p'_E}$ are tunable hyperparameters, $\mathbf{u}^V_{j,k}$, $\widehat{\mathbf{u}}^V_{j,k}$ are vectors of FOM outputs and GNN predictions for the output feature of the $j$-th node and for the $k$-th graph. Similarly, $\mathbf{u}^E_{i,k}$, $\widehat{\mathbf{u}}^E_{i,k}$ 
are vectors of FOM outputs and GNN predictions for the $i$-th edge output feature and for the $k$-th graph.
In particular, the GNN predictions $\widehat{\mathbf{u}}^V_{j,k}$ and $\widehat{\mathbf{u}}^E_{i,k}$ are feedforward mappings explicitly depending on the parameters $\theta$ that encode their trainable weights and biases.

The physics-informed component of the loss enforces the governing principles of the graph-based problem by penalizing the residuals of key equations~\cite{Raissi2019}:
\begin{equation}
    \mathcal{L}_{\text{physics}}(\theta) = \gamma \, \mathcal{L}_{\text{constitutive}}(\theta) + (1 - \gamma) \, \mathcal{L}_{\text{mass}}(\theta),
    \label{eq:physics_loss}
\end{equation}
where $\gamma \in [0,1]$ controls the balance between the constitutive relation and the mass conservation constraint.

The constitutive residual enforces the consistency between the predicted solution and the constitutive law governing the problem on the graph through the term:
\begin{equation}
    \mathcal{L}_{\text{constitutive}}(\theta) = 
    \frac{1}{N} \sum_{k=1}^N 
    \big\| C_k \widehat{\mathbf{U}}_k^V(\theta) - R^E_k \widehat{\mathbf{U}}_k^E(\theta) - \mathbf{F}^V_k \big\|^2.
    \label{eq:constitutive_loss}
\end{equation}
In the context of microvascular flow, this residual enforces, for example, the Poiseuille constitutive law, ensuring consistency between predicted pressure drops and flow rates along the vessels through the edge resistances computed from predicted features (e.g., diameter, hematocrit).

The mass balance residual enforces the continuity conditions at the junction nodes of the graph, penalizing imbalances of the nodal constraints:
\begin{equation}
    \mathcal{L}_{\text{mass}}(\theta) = 
    \frac{1}{N} \sum_{k=1}^N 
    \big\| C^{\!\top}_k \widehat{\mathbf{U}}^E_k(\theta) - \mathbf{F}^E_k \big\|^2.
    \label{eq:mass_loss}
\end{equation}
In the specific case of microvascular flow, this residual reduces to a mass conservation condition that ensures that total inflow and outflow at each bifurcation node are balanced.


\section{Microvascular blood flow models}
\label{sec:models}
Accurate modeling of blood flow in microvascular networks requires a careful combination of three key components: biologically meaningful vascular geometries (addressed in Section \ref{sec:networks}), physically consistent hemodynamic models, and computational methods capable of resolving flow across complex networks with high fidelity \cite{Popel2005,Secomb2017a}. In this section, we describe the modeling framework adopted to simulate blood flow through synthetic and anatomically realistic networks, formulated within the abstract setting introduced in Section~\ref{sec:abstract_metric_graph_problem}.
These simulations serve as the basis for the development and evaluation of the surrogate models presented in Section \ref{sec:method}. In particular, they are used both to generate training data for Graph Neural Network (GNN) models and to evaluate their predictive performance on out-of-distribution vascular topologies. 

\subsection{Physical models for microvascular blood flow}

The cornerstone of microvascular blood flow models is Poiseuille’s law, which describes the relationship between flow rate, pressure gradient, and hydraulic resistance in laminar flow conditions~\cite{Popel2005}.
Given a vessel of length $L$ and diameter $D$, the axial velocity of blood $v$ is related to the volumetric flow rate $Q$ by the expression:
\begin{equation}
    \label{eq:velocity_from_flow}
    v = \frac{Q}{\pi \frac{D^2}{4}} = \frac{\Delta P}{\pi R \frac{D^2}{4}},
\end{equation}
where $\Delta P$ denotes the pressure drop across the vessel and $R$ is the hydraulic resistance. According to Poiseuille’s law, resistance depends on the geometry of the vessel and the viscosity of the blood $\mu$ as follows:
\[
    R = \frac{128 \mu L}{\pi D^4}.
\]
For a given edge $e\in E$, this expression defines the coefficient $R_i^E$ and the constitutive function $r(\weight_i)$ in \eqref{eq.general_poiseuille}, being $\weight_i=\{\mu_i,D_i,L_i\}$.
Blood viscosity plays a central role in determining the resistance to flow and may be treated as either constant or variable depending on the physical regime considered. In microvascular domains, where RBC dynamics are non-negligible, the effective viscosity becomes a nonlinear function of both vessel diameter and local hematocrit (i.e. the volumetric concentration of RBCs in the blood \cite{Secomb2017a}).
Before detailing the main models, we list the key assumptions underlying our formulation \cite{GL,GL2,Popel2005,Secomb2017a}:
\begin{description}
\item[\textit{(i)}]  Blood flow and hematocrit distribution is stationary; 
\item[\textit{(ii)}]  RBC transport is dominated by advection;
\item[\textit{(iii)}]  Hematocrit is assumed to be uniform along the length of each vessel (no axial variations due to RBC production, loss, or plasma leakage);
\item[\textit{(iv)}]  Only Y-junctions with straight branches are considered, and angle dependence is neglected;
\item[\textit{(v)}]  The Fahraeus–Lindqvist and Zweifach–Fung  effects are explicitly included to capture key features of microvascular blood flow \cite{GL,GL2}, while other aspects such as leukocyte interactions and diameter irregularities are neglected.
\end{description}

We investigate two distinct modeling frameworks for microvascular blood flow. The first is a linear rheology model, which assumes a constant blood viscosity and neglects the effects of RBCs. This simplification leads to a linear system of equations such as \eqref{eq.general_poiseuille} that can be solved efficiently using standard numerical techniques. The second approach is a nonlinear rheology model, which incorporates the dependence of viscosity on hematocrit and accounts for the redistribution of RBC in vascular bifurcations. 
This added complexity results in a nonlinear system that requires iterative procedures to achieve convergence.

\subsubsection{Linear models}
The vascular domain is represented as a weighted directed graph
$\mathcal{G} = (V, E, W)$ defined before, where the weights on the edges $\weight^E_i=\{\mu_i,D_i,L_i\}$ are the diameter $D_i$ and the length $L_i$ of each edge and a uniform blood viscosity $\mu_i=\mu$. 
The hydraulic state of the network is characterized by the nodal vector $\mathbf{U}^V$ representing the pressure, that we name from now on $\mathbf{P} = \{P_j\}_{j=1}^{n}$ and an edge vector $\mathbf{U}^E$ representing the flow rates, named $\mathbf{Q} = \{Q_i\}_{i=1}^{m}$. Assuming laminar flow of a Newtonian fluid and neglecting red blood cells, the pressure drop along each edge follows a Poiseuille-like relation, while the nodal balance of flow enforces mass conservation:
\begin{equation}
    \label{eq:linear_solver}
    \left\{
    \begin{aligned}
    P_{s(i)} - P_{t(i)} &= R_i^E Q_i, & \quad i &= 1,\ldots,m,\\[2mm]
    \sum_{i\in E(v_j)} C_{ij} Q_i &= 0, 
    & \quad j &= 1,\ldots,n,
    \end{aligned}
    \right.
\end{equation}
The resistance of each vessel is given by $R_i^E = \tfrac{128\,\mu\,L_i}{\pi D_i^4}$, with $\mu = 3\;\mathrm{cP}$ denoting the dynamic viscosity of plasma at physiological temperature \cite{Popel2005}. 
Equations~\eqref{eq:linear_solver} describe a linear system reported in \eqref{eq:discrete_system} that defines the nodal pressures and edge flows on the graph $\mathcal{G}$. 
The velocities in the vessels can then be computed using Equation~\eqref{eq:velocity_from_flow}. 
To manage the boundary conditions, we first label the vertices $V$ as interior, inlet and outlet nodes.
We set $V_{inlet}:=\{v_j\}_{j \in J_{inlet}}$ as $J_{inlet}$ is a multi-index identifying the input nodes. We proceed in a similar way for the set of outlet nodes $V_{outlet}$.
As boundary values, we prescribe Dirichlet conditions, assigning the pressures at the inlet nodes and the outlet nodes named $\mathbf{P}^{in}:=\{P^{in}_j\}_{j\in J_{inlet}}$ and $\mathbf{P}^{out}:=\{P^{out}_j\}_{j\in J_{outlet}}$, respectively.

We can now assemble a nodal-valued boundary values vector $\overline{\mathbf{P}}$, identifying the weights on the nodes $\weight_j^V=\overline{P}_j$, and a diagonal matrices $B^{in},B^{out}\in\mathbb{Z}^{n\times n}$, encoding the inlet and outlet Dirichlet boundary nodes, respectively, such that:

\[ \overline{P}_j = 
    \left\{ \begin{array}{cl}
            P^{in}_j, & \forall j \in J_{inlet} \\
            P^{out}_j, & \forall j \in J_{outlet}  \\
            0, & \text{otherwise}
    \end{array} \right.
    \quad
    B_{jj}^{in}:=
    \left\{ \begin{array}{cl}
            1,  & \forall j \in J^{inlet} \\
            0 & \text{otherwise}
    \end{array} \right.
   \quad
    B_{jj}^{out}:=
    \left\{ \begin{array}{cl}
            1,  & \forall j \in J^{outlet} \\
            0 & \text{otherwise}
    \end{array} \right.
\]

This entails that we can derive the following algebraic formulation for \eqref{eq:linear_solver}:
\begin{equation}\label{eq:discrete_system_lin}
    L^W \, \mathbf{U}^W = \mathbf{F}^W,
    \quad 
    L^W :=
    \begin{bmatrix}
        C & -R^E \\[10pt]
        B^{in}+B^{out} & (I-B^{in}-B^{out})C^{\!\top}
    \end{bmatrix},
    \quad
    \mathbf{F}^W:=
    \begin{bmatrix}
        \mathbf{0} \\[10pt]
        \overline{\mathbf{P}} 
    \end{bmatrix},
\end{equation}

where the system in \eqref{eq:discrete_system} is assembled with constant viscosity, leading to a fixed resistance matrix $R^E$, and is solved using a direct sparse linear solver. In our implementation, the linear system is assembled as a sparse matrix, converted to compressed sparse row (CSR) format, and solved on the CPU through a single call to \texttt{scipy.sparse.linalg.spsolve}.
Within the linear rheology framework, the solution operator is denoted by $\mathcal{S}_L$, such that $\mathcal{S}_L(G, W) = \mathbf{U}^W$.

\subsubsection{Nonlinear models}
\label{sec:nonlinear_models}

Microvascular blood flow exhibits nonlinear rheological behavior due to two main factors: the dependence of blood viscosity on hematocrit and vessel diameter (known as the Fåhræus–Lindqvist effect) and the heterogeneous distribution of RBCs in bifurcations, commonly referred to as plasma skimming (also called the Zweifach-Fung effect) \cite{GL,GL2}. To capture these effects, we adopt a biphasic nonlinear flow model in which blood is treated as a suspension of plasma and RBCs. 
Recalling that the hematocrit is the volumetric concentration of RBCs in the blood, we augment the unknowns at the edges including the hematocrit levels. As a result, in this model the unknowns are the pressures on the nodes $\mathbf{U}^V:=\mathbf{P}=\{P_j\}_{j=1}^n$, the flow rates, and the hematocrit on the edges, i.e. $\mathbf{Q}=\{Q_i\}_{i=1}^m,\,\mathbf{H}=\{H_i\}_{i=1}^m$, respectively, so that $\mathbf{U}^E:=\{ [Q_i, H_i] \}_{i=1}^m$.
In the nonlinear rheology regime, the governing equations define a solution operator $\mathcal{S}_N(G, W) = \mathbf{U}^W$, which maps the graph and the physical parameters to pressure, flow rate and hematocrit.

Similarly to the linear setting, we introduce the diagonal selector matrices $B^{in}, B^{out} \in \mathbb{Z}^{n\times n}$ to handle boundary nodes. The corresponding boundary vectors $\overline{\mathbf{P}}$ and $\overline{\mathbf{H}}$ collect the inlet and outlet pressure and inlet hematocrit values, respectively.

The effective viscosity $\mu(D,H)$ varies with the local discharge hematocrit $H$ and the diameter of the vessel $D$, incorporating the Fåhræus–Lindqvist effect. Moreover, the hematocrit field evolves under the influence of a kinematic plasma skimming model, ensuring consistent RBC mass conservation across bifurcations \cite{GL,GL2,linningerplos,VL1}.
The coupled nonlinear system is solved through a fixed-point iterative scheme alternating between two subproblems: (i) a pressure-flow problem with hematocrit-dependent viscosity and (ii) a convection-like problem for RBC transport governed by pseudo-fluxes. This algorithm is illustrated in the following.

\paragraph{Blood flow with hematocrit-dependent viscosity}
Given an estimate of the hematocrit values at each edge of the graph $\mathbf{H}=\{H_i\}_{i=1}^m$ and the weights on the vertices $\weight^V_j=\{\overline{P}_j, \overline{H}_j\}$ and on the edges $\weight^E_i=\{D_i, L_i\}$, we assemble the problem \eqref{eq:discrete_system} with 
\begin{equation*}
    R_i^E = r(\weight_i,L_i) =\frac{128\, \mu_i(D_i, H_i)\, L_i}{\pi D_i^4},
\end{equation*}
where resistances are calculated according to the empirical law of Pries et al.~\cite{Pries1992},
\begin{equation}
    \mu_i(D_i, H_i) = \mu_{\text{plasma}} \left[1 + (\mu_{0.45,i} - 1) \frac{(1 - H_i)^{\gamma_i} - 1}{(1 - 0.45)^{\gamma_i} - 1} \right],
    \label{eq:pries_mu}
\end{equation}
where $\mu_{\text{plasma}} = 1$ cP is the plasma viscosity, $L_i$ is the length of vessel $i$, and the parameters $\mu_{0.45,i}$ and $\gamma_i$ depend on $D_i$ as:
\begin{align*}
    \mu_{0.45,i} &= 220 e^{-1.3 D_i} + 3.2 - 2.44 e^{-0.06 D_i^{0.645}}, \\
    \gamma_i &= \left(0.8 + e^{-0.075 D_i} \right) \Big(-1 + \frac{1}{1 + 10^{-11} D_i^{12}}  \Big) + \frac{1}{1 + 10^{-11} D_i^{12}}.
\end{align*}
This problem is solved using the same matrix formulation as in the linear case, namely \eqref{eq:discrete_system}, with hematocrit-dependent resistances.

\paragraph{Uneven hematocrit and convection formulation} In the second stage of this scheme, we update the hematocrit distribution throughout the network while preserving the mass of RBC. The conservation of RBCs in a bifurcation involving one parent vessel $(0)$ and two daughter vessels $(1)$ and $(2)$  can be expressed as follows (these are local indices, not to be confused with the global numbering of the graph $\mathcal{G}$):
\begin{equation}
    Q_{(0)} H_{(0)} = Q_{(1)} H_{(1)} + Q_{(2)} H_{(2)}.
    \label{eq:rbc_mass_balance}
\end{equation}

However, due to the strong nonlinear coupling between $Q$ and $H$—especially through the hematocrit-dependent viscosity law introduced in Eq.~\eqref{eq:pries_mu}—direct application of Eq.~\eqref{eq:rbc_mass_balance} across large bifurcating microvascular networks can lead to numerical instability or nonphysical solutions. To mitigate this, we introduce a virtual node potential $H^*$ at each bifurcation node~\cite{GL}. This reformulation allows plasma skimming to be posed as a linear convection problem over $\mathcal{G}$, which ensures the uniqueness and stability of the solution \cite{ventimiglia2023meshfree}. 

Let us formulate first the model locally, i.e. at the level of one bifurcation, and then extend it to the graph level.
The discharge hematocrit $H_{(i)}$ in each downstream daughter branch $i=1,2$ is then determined by:
\begin{equation}
    H_{(i)} = \lambda_{(i)} H^*,
    \label{eq:Hface}
\end{equation}
where $\lambda_{(i)}$ is a kinematic plasma skimming coefficient. This coefficient captures the tendency of RBCs to preferentially enter one daughter branch over another, and is defined as a function of the ratio of daughter-to-parent diameters
$\lambda_{(i)} = \left( {D_{(i)}}/{D_{(0)}} \right)^{2/M_d}$,
where $D_{(i)}$ and $D_{(0)}$ are the diameters of the daughter and parent vessels, respectively, and $M_d$ is the drift parameter. In this work, we set $M_d = 5.25$ based on physiological evidence.

Substituting Eq.~\eqref{eq:Hface} into Eq.~\eqref{eq:rbc_mass_balance} yields the following,
\begin{equation*}
    Q_{(0)} H_{(0)} = Q_{(1)} \lambda_{(1)} H^* + Q_{(2)} \lambda_{(2)} H^* = Q_{(1)}^* H^* + Q_{(2)}^* H^*,
    \label{eq:rbc_balance_virtual}
\end{equation*}
where the geometry-adjusted pseudo-fluxes $Q_{(i)}^*$ are defined as $Q_{(i)}^* = Q_{(i)} \lambda_{(i)}$.
Solving for $H^*$ gives:
\begin{equation}
    H^* = \frac{Q_{(0)}^*}{Q_{(1)}^* + Q_{(2)}^*} \overline{H}^*,
    \label{eq:Hstar}
\end{equation}
where $\overline{H}^*$ is the hematocrit potential of the upstream node from which the parent vessel originates.

Equation~\eqref{eq:Hstar} enables a topologically ordered computation of node hematocrit potentials across the network. 
At inflow nodes, we prescribe $H^* = H^{in}$ and set $\lambda_i = 1$ for all incident edges. Once all $H^*$ values are computed, the discharge hematocrits $H_i$ along vessels are recovered via Eq.~\eqref{eq:Hface}.

Let us now discuss how to formulate the previous algorithm at the global level, namely for the entire graph $\mathcal{G}$.
The discrete convection equation for the vector of hematocrit node potentials $H^* \in \mathbb{R}^n$ can be expressed in matrix form:
\begin{equation}
    N(\mathbf{Q}^*) \mathbf{H}^* = N^{in}(\mathbf{Q}^*) \overline{\mathbf{H}}.
    \label{eq:convection_matrix}
\end{equation}
The matrices $N(\mathbf{Q}^*), N^{in}(\mathbf{Q}^*) \in \mathbb{R}^{n \times n}$ are defined as:
\begin{align*}
    N(\mathbf{Q}^*) &= C^\top \cdot \mathrm{diag}(\mathbf{Q}^*) \llbracket -C \rrbracket - \mathrm{diag}(B^{out} C^\top \mathbf{Q}^*), \\
    N^{in}(\mathbf{Q}^*) &= - \mathrm{diag}(B^{in} C^\top \mathbf{Q}^*),
\end{align*}
where $C$ is the incidence matrix of the graph and $\llbracket X \rrbracket_{ij} = \max(X_{ij}, 0)$ denotes the positive-part operator. 

In vascular confluences (e.g., venular unions), we set $\lambda = 1$ by definition. A more detailed explanation of this discretization can be found in~\cite{VL1} and~\cite{ventimiglia2023meshfree}, where it is derived from a finite volume formulation of convection in graphs.
Once node potentials are determined, the hematocrit values along the vessel segments are obtained by solving the following system,
\begin{equation}
    \mathbf{H} = \Lambda \llbracket -C \rrbracket \mathbf{H}^*,
    \label{eq:recover_H}
\end{equation}
where $\mathbf{H} \in \mathbb{R}^m$ is the vector of hematocrits defined on edges and $\Lambda \in \mathbb{R}^{m \times m}$ is a diagonal matrix with the $\lambda_i$ coefficients for each edge.

\medskip

The full nonlinear model is then solved via fixed-point iterations, alternating between the following steps:
\begin{description}
    \item[(i)] computing the flow field $\mathbf{Q}$ based on hematocrit-dependent viscosity;
    \item[(ii)] updating node potentials $\mathbf{H}^*$ via Eq.~\eqref{eq:convection_matrix} and recovering edge values $\mathbf{H}$ from Eq.~\eqref{eq:recover_H}.
\end{description}

The iteration is initialized with a uniform hematocrit field (e.g., $H_i = 0.4$) and repeated until both flow and hematocrit fields converge. At each fixed-point iteration, the corresponding linear systems are solved on the CPU using the same sparse direct solver employed for the linear problem.
Convergence is typically achieved within 5–7 iterations, even for large networks. A key computational advantage of this approach is that the convection matrix $N(\mathbf{Q}^*)$ is lower triangular when the nodes are ordered topologically with respect to the flow direction. This guarantees that Eq.~\eqref{eq:convection_matrix} has a unique solution, provided that no edge flows vanish. The graph-based convection framework significantly improves numerical stability compared to classical models based on recursive split rules \cite{Pries1992,pries2005microvascular}.

The linear and nonlinear models described above provide a rigorous computational framework to simulate microvascular blood flow across a wide range of regimes. However, their direct application to large ensembles of vascular graphs is computationally demanding, especially in the nonlinear setting where iterative convergence is required.


\section{Mathematical synthesis of vascular networks}
\label{sec:networks}
The construction of surrogate models for microcirculatory flow requires access to physiologically meaningful vascular graphs that can be efficiently handled in a digital setting. In this Section, we focus on the generation of such networks, presenting two complementary approaches that differ in terms of anatomical fidelity, computational cost, and scalability.  
Section~\ref{sec:synthetic_general} introduces a fast geometry-based procedure that produces large collections of simplified but topologically rich capillary networks, suitable for training data-driven surrogates.  
Section~\ref{sec:cco}, instead, describes an optimization-based strategy for image-based Cerebral Network Synthesis, which produces high-fidelity cerebrovascular architectures that serve as rigorous benchmarks for validation.  
These two approaches synergistically interact with the generation of GNNs-based surrogate models of blood flow, as we train the GNNs on light vascular graphs obtained through the approach of Section~\ref{sec:synthetic_general}, while we test them on the anatomically accurate configurations of Section~\ref{sec:cco}.

%

\subsection{Fast generation of general-purpose synthetic capillary networks}\label{sec:synthetic_general}
The core of this algorithm uses Voronoi 3D tessellations and the Dijkstra algorithm to obtain a graph $G^s\in\mathbb{Z}^{m \times n}$ that represents a viable vascular network. The procedure for generating the synthetic networks can be summarized in the following steps: (i) generation of the vascular network in a cubic domain; (ii) removal of all trifurcations; (iii) assignment of diameters; (iv) rescaling the cube dimension to have a coherent Surface-to-Volume (S/V) ratio.

\paragraph{Vascular network generation.} Synthetic vascular networks are generated within a unit cubic domain, subdivided into three horizontal layers to introduce heterogeneity throughout depth. In each layer, a point cloud is sampled according to a predefined density function, controlling the spatial variation in vascular density. The aggregate point set is used as input for a three-dimensional Voronoi tessellation, whose edges define the centerlines of the candidate vessel. Unfeasible branches, e.g., isolated or highly tortuous segments, are subsequently removed through filtering procedures. Inlet and outlet planes are then identified geometrically; points in proximity to these planes are projected onto them and connected via straight segments, establishing well-defined boundary inflow and outflow conditions. The resulting vascular graph exhibits structural complexity and spatial variability, while preserving computational feasibility. The number of outlets is about twice as many as the number of inlets.
    
\paragraph{Removal of trifurcations.} 
Capillary networks can commonly be considered isotropic
anastomosing networks, whose vertices typically have up to three connections~\cite{Duvernoy,Blinder2013}.
To enforce a tree-like topology and eliminate trifurcations, a shortest-path strategy is employed via Dijkstra’s algorithm. Starting from selected inlet nodes, the algorithm traces minimal paths to designated outlet nodes, assembling a spanning structure that avoids cycles and guarantees flow feasibility. Detected trifurcations are resolved by removing the shortest branch, inserting an auxiliary node on the adjacent vessel, and reconnecting the remaining segments to preserve graph connectivity.

\paragraph{Assignment of diameters.} Vessel diameters are assigned following a multi-stage algorithm. Initially, diameters are drawn randomly, following a uniform distribution, to enable a preliminary solution of the linear flow model, which provides reliable orientation for topological sorting. The maximum diameter is then assigned to the input segments and propagated downstream with progressive attenuation at the bifurcation points. For Y-shaped bifurcations, a modified version of Murray’s law is adopted, namely:
$r^4 = \sqrt[4]{r_1^4 + r_2^4}$
as proposed in \cite{akita2016experimental}, to better represent microvascular scaling laws and mitigate overproduction of vessels of minimal diameter. 
Throughout, a minimum and maximum diameter threshold is enforced. 
If $r_1$ and $r_2$ were known, $r$ is found applying the formula. On the other hand, if only $r$ is known, $r_1$ is randomly generated and then $r_2$ is calculated. In each computation, a check that the diameter is between a minimum and a maximum value is performed. 

\paragraph{Rescaling the cube dimension.} 
As a final step, the domain size is rescaled to ensure a physiologically realistic surface-to-volume ($S/V$) ratio. Specifically, given the total vascular surface $S = \sum_{i \in E} \pi D_i L_i$ the cube volume $V$ is adjusted such that $S/V \approx 7000, \mathrm{m^{-1}}$, in alignment with physiological observations. Figure~\ref{fig:arterialtree} illustrates a representative output of the synthetic network generator.

\subsection{Anatomically accurate cerebral vascular networks}
\label{sec:cco}

To evaluate the transfer of the proposed GNN pressure predictions, we test its performance on anatomically accurate vascular networks generated through a constrained optimization framework, named the image-based Cerebral Network Synthesis (iCNS)~\cite{linninger1, linninger2}, and encoded by $G^a\in\mathbb{Z}^{m \times n}$. These networks are constructed using an image-based extension of the Constrained Constructive Optimization (CCO) methodology, originally proposed for coronary arterial tree synthesis~\cite{SCHREINER199927,KARCH199919}.

The CCO paradigm incrementally grows a vascular tree by connecting newly sampled terminal nodes to existing segments of the network, solving at each iteration a constrained minimization problem. The objective is to minimize the total volume of the tree as a function of the position of the bifurcation $x \mapsto V(x)$, while enforcing global hemodynamic consistency in terms of continuity of pressure and flow conservation between bifurcations.
Formally, the problem can be cast as a nonlinear optimization problem involving: 
\begin{description}
    \item[\textit{(i)}] \textit{topological decisions}: the choice of the parent segment to which a new terminal node is attached;
    \item[\textit{(ii)}] \textit{geometric optimization}: the placement of the bifurcation point determining local segment lengths and diameters.
\end{description}

\begin{figure}[t]
    \centering
    \includegraphics[width=1\textwidth]{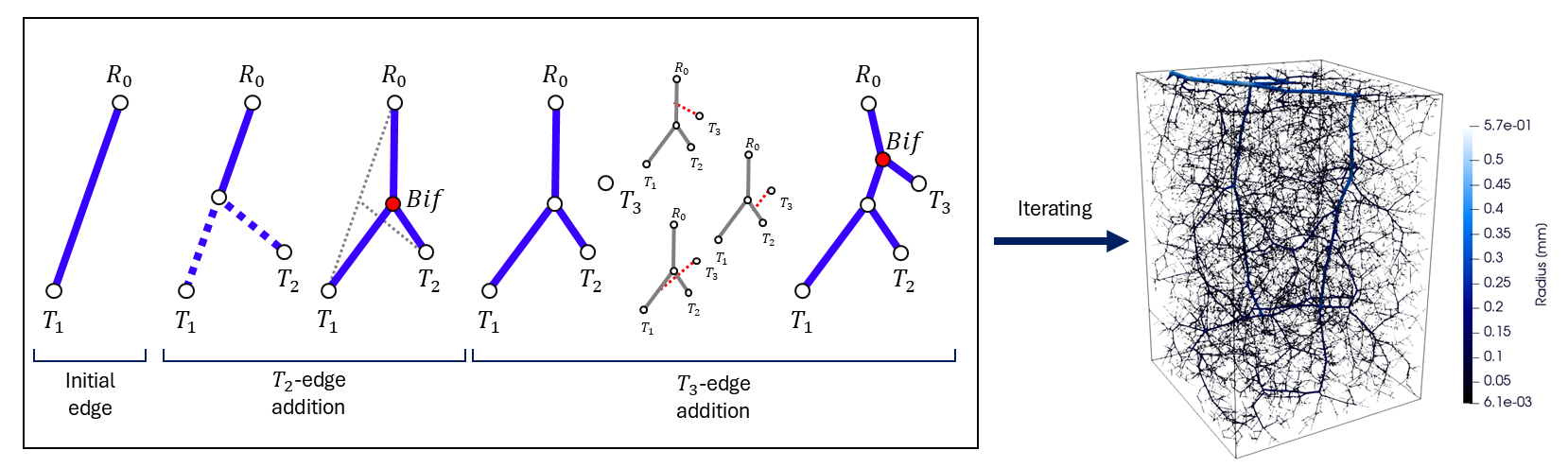}
    \caption{Example of an anatomically guided synthetic arterial tree generated via iCNS, showing the bifurcation optimization process in the CCO framework for each segment addition, adapted from \cite{linninger1, linninger2}.}
    \label{fig:arterialtree}
\end{figure}

Such a problem is intrinsically NP-hard due to the coupling between discrete (connectivity) and continuous (geometry) variables. We refer the reader to Figure \ref{fig:arterialtree} for a visualization of these optimization steps.
To ensure anatomical fidelity, the iCNS method is constrained within volumetric regions derived from imaging data, such as micro-CT or mouse brain atlases. These provide anatomical priors that accommodate anatomical realism.
Through recursive optimization, the algorithm produces large-scale vascular networks that incorporate arterial, capillary, and venous compartments. A final stage of microvascular closure connects the terminal arterial and venous nodes with tortuous capillary segments, resulting in a fully physiologically connected and stable structure. The resulting graphs reproduce realistic morphometrics, including vessel density, length distribution, diameter hierarchy, and topological depth, and are suitable for the validation of microcirculatory flow models (Figure \ref{fig:arterialtree}).

\section{Surrogate blood flow models based on GNNs}
\label{sec:gnn_models}
We introduce here four GNN configurations designed to approximate the parameter-to-solution map 
$\mathcal{S}(G, W) = \mathbf{U}^W = [\mathbf{U}^V, \mathbf{U}^E]^{\top}$, for the blood flow models introduced in Section~\ref{sec:models}.
Each configuration differs in the complexity of its physical modeling and the composition of its loss function, but all share the same set of input features defined on the nodes and edges of the graph. Table \ref{tab:parameters} summarizes the features of the four GNN models that are described in what follows.

The GNN models are trained on ensembles of $N$ synthetic vascular graphs $G^s_k = (V_k, E_k), \ k=1,...,N$, generated using the geometry-based procedure described in Section~\ref{sec:synthetic_general}. Since these vascular graphs can be created at a low cost and the blood flow models can be efficiently resolved on these setups, this option simplifies the implementation of the supervised learning problem for training the surrogate GNN models.

The input feature vector $\widetilde{W} = [\widetilde{W}^V, \widetilde{W}^E]^{\top}$ contains geometric and physical quantities. 
Specifically, the node features vector $\widetilde{\weight}_j^V$ includes boundary pressures ($P_j^{\text{in}}, P_j^{\text{out}}$) at inlet and outlet vertices together with their extrema 
($\max(P^{\text{in}})$, $\min(P^{\text{out}})$).
The edge features $\widetilde{W}^E$ include the vessel length $L_i$ and diameter $D_i$. 
For the nonlinear rheology model, the hematocrit boundary value $H^{\text{in}}$ is also included among the edge inputs. 
As a result, we have
\[
\widetilde{\weight}_j^V=\{P_j^{\text{in}}, P_j^{\text{out}},\max(\mathbf{P}^{\text{in}}), \min(\mathbf{P}^{\text{out}})\},
\quad
\widetilde{\weight}_i^E=\{D_i,L_i,H^{\text{in}}\}.
\]

The output feature vector $\mathbf{U}^W = [\mathbf{U}^V, \mathbf{U}^E]$ consists of the predicted nodal and edge quantities, which, depending on the model, include pressure $P_j$, transformed velocity $v_i = \mathcal{T}_v(Q_i)$, and hematocrit $H_i$. More precisely, from the calculated flow rates, we evaluate the corresponding blood velocities along each segment of the vessel. To mitigate the skewness of the velocity distribution and improve the numerical stability of the learning process, we introduce a \emph{logarithmic velocity transformation}
$\mathcal{T}_v : \mathbb{R} \to \mathbb{R}$, defined as:
\[
v:=\mathcal{T}_v(Q)
= \frac{\operatorname{sign}(Q)}{k_v} 
\log\!\left(1 + \frac{4 |Q|}{\pi D^2}\right),
\quad
\mathcal{T}_v^{-1}(v) =
\operatorname{sign}(v)\;
\frac{\pi D^2}{4} 
\left( e^{k_v |v|} - 1 \right),
\]
where $k_v$ is an empirically chosen scaling constant. The inverse mapping $\mathcal{T}_v^{-1} : \mathbb{R} \to \mathbb{R}$ 
is used to recover the physical flow rate values from the transformed representation.
Unless explicitly stated otherwise, the velocity errors reported below
are computed after mapping both the GNN predictions and the reference
outputs back to physical flow rate and axial velocity, rather than on the transformed variables used during training.
Comprehensively, the output features of the nonlinear rheology model are the following:
\[
\mathbf{u}_j^V=\{P_j\}, 
\quad
\mathbf{u}_i^E=\{v_i,H_i\},
\]
while in the linear models the hematocrit term $H_i$ is not included among the outputs.

\begin{table}[t]
\centering
\renewcommand{\arraystretch}{1.2}
\begin{tabular}{|c|c|c|c|}
\hline
\multicolumn{4}{|c|}{\textbf{Input Features}} \\
\hline
\textbf{Symbol} & \textbf{Description} & \textbf{Type} & \textbf{Models} \\
\hline
$L_i$ & Vessel length & Edge feature & 1--4 \\
$D_i$ & Vessel diameter & Edge feature & 1--4 \\
$P_j^{in}$, $P_j^{out}$ & Inlet/outlet pressures & Node feature & 1--4 \\
$\max(\mathbf{P}^{in})$, $\min(\mathbf{P}^{out})$ & Maximum/Minimum boundary pressure & Node feature & 1--4 \\
$H^{in}$ & Hematocrit boundary value & Edge feature & 4 \\
\hline
\multicolumn{4}{|c|}{\textbf{Output Features}} \\
\hline
$P_j$ & Pressure & Node feature & 1--4 \\
$v_i$ & Transformed velocity & Edge feature & 2--4 \\
$H_i$ & Hematocrit & Edge feature & 4 \\
\hline
\end{tabular}
\caption{Summary of input and output features for the GNN models.}
\label{tab:parameters}
\end{table}

\subsection{GNN model variants}
The configuration of the four GNN-based surrogate models is described in the following, with increasing complexity from the simplest to the most complex and physically accurate model.

\paragraph{\textbf{Model 1:} Vanilla GNN (Linear Rheology, Pressure Only)}
This model approximates only the nodal pressure field $\widehat{\mathbf{U}}^V(\theta^*)=\widehat{\mathbf{P}}_k(\theta^*)=\widehat{\mathcal{S}}_1(G^s_k,\widetilde{W}_k;\theta^*)\approx \mathcal{S}_L(G^s_k,W_k)$, with $\widehat{\mathbf{U}}^V\in \mathbb{R}^{n_k}, \ \widehat{\mathbf{U}}^E=\emptyset$, 
on each graph $G^s_k$, using features: 
\[
\widetilde{\weight}_j^V=\{P_j^{\text{in}}, P_j^{\text{out}},\max(\mathbf{P}^{in}), \min(\mathbf{P}^{out})\},\, j=1,\ldots,n_k,
\quad
\widetilde{\weight}_i^E=\{D_i,L_i\},\,\, i=1,\ldots,m_k.
\]
Given labeled data $\mathbf{P}_k=\mathcal{S}_1(G^s_k,{W}_k) \in \mathbb{R}^{n_k}$
we use a mean-square supervised loss:
\[
\mathcal{L}(\widehat{\mathcal{S}}_1,\mathbf{U}^W; \theta) = \mathcal{L}_{\text{pressure}}(\theta)
= \frac{1}{N} \sum_{k=1}^N \frac{1}{n_k} 
\left\| \mathbf{P}_k - \widehat{\mathbf{P}}_k(\theta) \right\|^2. 
\]
The optimal parameters $\theta$ that describe the FFNN in every module of the GNN architecture are determined by solving the following minimization problem (the total number of parameters for this configuration, as well as for the following cases, is reported in Table \ref{tab:model_hyperparams_sensitivity}):
\[
\theta^*=\argminF_{\theta \in \Theta} \mathcal{L}(\widehat{\mathcal{S}}_1,\mathbf{U}^W; \theta).
\]

\paragraph{\textbf{Model 2:} Data-Driven GNN (Linear Rheology, Pressure + Velocity)}
This configuration extends the previous one by jointly predicting nodal pressures and edge velocities
\[
[\widehat{\mathbf{U}}^V(\theta),\widehat{\mathbf{U}}^E(\theta)]^{\top}=[\widehat{\mathbf{P}}_k(\theta),\widehat{\mathbf{v}}_k(\theta)]^\top=\widehat{\mathcal{S}}_2(G^s_k,\widetilde{W}_k;\theta)\approx \mathcal{S}_L(G^s_k,W_k),
\]
using the same input features of Model 1.
Given labeled data on nodes and edges
$[\mathbf{P}_k,\mathbf{v}_k]^{\top}=\mathcal{S}(G_k,{W}_k) \in \mathbb{R}^{n_k+m_k}$, we retrieve the minimizer $\theta^*$ of the loss that combines the two supervised terms:
\[
\mathcal{L}(\widehat{\mathcal{S}}_2,\mathbf{U}^W; \theta) =  \gamma_D \mathcal{L}_{\text{pressure}}(\theta) + (1 - \gamma_D) \mathcal{L}_{\text{velocity}}(\theta), \quad \gamma_D \in [0,1],
\]
\[
\mathcal{L}_{\text{pressure}}(\theta) = \frac{1}{N} \sum_{k=1}^N \frac{1}{n_k} \left\| \mathbf{P}_k - \widehat{\mathbf{P}}_k(\theta) \right\|^2,
\qquad 
\mathcal{L}_{\text{velocity}}(\theta) = \frac{1}{N} \sum_{k=1}^N \frac{1}{m_k} \left\| \mathbf{v}_k - \widehat{\mathbf{v}}_k(\theta) \right\|^2.
\]

\paragraph{\textbf{Model 3:} Physics-Informed GNN (Linear Rheology)}
In this model, we incorporate physics-based regularization to enforce the constitutive and conservation relations of the flow problem.
The GNN predicts both nodal pressures and edge velocities in all nodes and edges of the vascular graph,
\[
[\widehat{\mathbf{U}}^V(\theta),\widehat{\mathbf{U}}^E(\theta)]^{\top}
=[\widehat{\mathbf{P}}_k(\theta),\widehat{\mathbf{v}}_k(\theta)]^{\top}
=\widehat{\mathcal{S}}_3(G^s_k,\widetilde{W}_k;\theta)
\approx \mathcal{S}_L(G^s_k,W_k),
\]
and is trained by minimizing a composite loss that combines data fidelity and physical consistency:
\[
\mathcal{L}(\widehat{\mathcal{S}}_3,\mathbf{U}^W;\theta) =
\delta \mathcal{L}_{\text{data}}(\theta) + (1-\delta)\mathcal{L}_{\text{physics}}(\theta),
\qquad \delta \in [0,1].
\]
The supervised term $\mathcal{L}_{\text{data}}$ is identical to Model~2, while the physics-informed component is defined as
\[
\mathcal{L}_{\text{physics}}(\theta) = \gamma_P \mathcal{L}_{\text{constitutive}}(\theta) + (1 - \gamma_P) \mathcal{L}_{\text{mass}}(\theta),
\qquad \gamma_P \in [0,1].
\]
The two physical residuals read:
\begin{align*}
\mathcal{L}_{\text{constitutive}}(\theta) 
&= \frac{1}{C_P} \frac{1}{N} 
\sum_{k=1}^{N} \frac{1}{m_k} 
\left\| 
C\,\widehat{\mathbf{P}}_k(\theta) - R^E \mathcal{T}_v^{-1}( \widehat{\mathbf{v}}_k(\theta))
\right\|^2, \\[4pt]
\mathcal{L}_{\text{mass}}(\theta) 
&= \frac{1}{C_M} \frac{1}{N} 
\sum_{k=1}^{N} \frac{1}{n_k}
\left\|
(I - B^{in} - B^{out}) C^{\!\top} \mathcal{T}_v^{-1}( \widehat{\mathbf{v}}_k(\theta))
\right\|^2,
\end{align*}
where $C_P = 35\,\text{mmHg}$ and $C_M = 10^6$ are normalization constants.
This formulation enforces the Poiseuille constitutive law at the edges of the vascular network and conservation of mass at the junctions of the network.

\paragraph{\textbf{Model 4:} Physics-Informed GNN (Nonlinear Rheology)}
This most comprehensive configuration extends the previous model by accounting for nonlinear viscosity and hematocrit-dependent effects.
The network outputs nodal pressures, edge velocities, and hematocrit levels,
\[
[\widehat{\mathbf{P}}_k(\theta), \widehat{\mathbf{v}}_k(\theta), \widehat{\mathbf{H}}_k(\theta)]^{\top}
=\widehat{\mathcal{S}}_4(G^s_k,\widetilde{W}_k;\theta)
\approx \mathcal{S}_N(G^s_k,W_k).
\]
The total loss function extends that of Model~3 by including hematocrit-dependent terms in both data-based and physics-based components.
The data-based contribution reads:
\[
\mathcal{L}_{\text{data}}(\theta) =
\frac{\gamma_{D_1} + \gamma_{D_2}}{2} \mathcal{L}_{\text{pressure}}(\theta) +
\frac{1 - \gamma_{D_1}}{2} \mathcal{L}_{\text{velocity}}(\theta) +
\frac{1 - \gamma_{D_2}}{2} \mathcal{L}_{\text{hematocrit}}(\theta),
\]
with $\gamma_{D_1}, \gamma_{D_2} \in [-1, 1],
\ \gamma_{D_1}+\gamma_{D_2} \geq 0,$ and where

\[
\mathcal{L}_{\text{hematocrit}}(\theta) = \frac{1}{N} \sum_{k=1}^N \frac{1}{m_k} \left\| \mathbf{H}_k - \widehat{\mathbf{H}}_k(\theta) \right\|^2.
\]
The physics-informed term includes the nonlinear resistance 
$R^E(\widehat{\mathbf{H}})$ and the hematocrit-weighted conservation residuals:

\[
\mathcal{L}_{\text{physics}} = \frac{\gamma_{P_1} + \gamma_{P_2}}{2} \mathcal{L}_{\text{constitutive}} + \frac{1 - \gamma_{P_1}}{2} \mathcal{L}_{\text{mass}}^{(1)} +
\frac{1 - \gamma_{P_2}}{2}
\mathcal{L}_{\text{mass}}^{(2)},
\]
with $\gamma_{P_1}, \gamma_{P_2} \in [-1, 1],
\ \gamma_{P_1}+\gamma_{P_2} \geq 0,$ and where
\begin{align*}
\mathcal{L}_{\text{constitutive}}(\theta) 
&= \frac{1}{C_P} \frac{1}{N} 
\sum_{k=1}^{N} \frac{1}{m_k} 
\left\| 
C\,\widehat{\mathbf{P}}_k(\theta) - R^E(\widehat{\mathbf{H}}_k(\theta)) \mathcal{T}_v^{-1}(\widehat{\mathbf{v}}_k(\theta))
\right\|^2, \\[4pt]
\mathcal{L}_{\text{mass}}^{(1)}(\theta) 
&= \frac{1}{C_M} \frac{1}{N} 
\sum_{k=1}^{N} \frac{1}{n_k}
\left\|
(I - B^{in} - B^{out}) C^{\!\top} \mathcal{T}_v^{-1}(\widehat{\mathbf{v}}_k(\theta))
\right\|^2, \\[4pt]
\mathcal{L}_{\text{mass}}^{(2)}(\theta) 
&= \frac{1}{C_M} \frac{1}{N} 
\sum_{k=1}^{N} \frac{1}{n_k}
\left\|
(I - B^{in} - B^{out}) C^{\!\top} 
\big(\widehat{\mathbf{H}}_k(\theta) \odot \mathcal{T}_v^{-1}(\widehat{\mathbf{v}}_k(\theta))\big)
\right\|^2.
\end{align*}
Here, the term $\mathcal{L}_{\text{mass}}^{(2)}$ ensures hematocrit-weighted flow conservation, while the nonlinear dependence $R^E(\widehat{\mathbf{H}})$ captures the Fåhræus–Lindqvist viscosity effect along each edge.

\subsection{Numerical setup}
Training and evaluation of GNN models are carried out using synthetic vascular networks generated by the algorithm detailed in Section~\ref{sec:models}. 
We generate a dataset of $N = 1200$ graphs with randomized topology and geometry, partitioned as follows: (i) \textit{Training set:} $N_\text{train} = 960$ graphs; (ii) \textit{Validation set:} $N_\text{val} = 120$ graphs; (iii) \textit{Test set:} $N_\text{test} = 120$ graphs.

Each network comprises approximately 300 nodes and 400 edges, with 25 to 40 inlets located at the top surface of the domain, corresponding to about 50–80 outlets. Boundary pressures are randomly sampled within physiological ranges, and all internal quantities are normalized to $[0, 1]$ to ensure numerical stability. 

\begin{table}[t]
\centering
\scriptsize
\renewcommand{\arraystretch}{1.12}
\setlength{\tabcolsep}{3.2pt}
\definecolor{lightgreenrow}{RGB}{217,234,211}

\begin{tabular}{|c|c|c|c|c|c|}
\hline
\textbf{Model} & \textbf{Rheology} & \textbf{\# Param.} & \textbf{Loss hyperparam.} & \textbf{$L^2$ pres. error} & \textbf{$L^2$ vel. error} \\
\hline
Model 1 & Linear & 31,144 & -- & 2.94\% & -- \\
\hline
Model 2 & Linear & 73,992 & $\delta=1,\ \gamma_D=0.4$ & 5.73\% & 36.7\% \\
\rowcolor{lightgreenrow}
Model 2 & Linear & 73,992 & $\delta=1,\ \gamma_D=0.5$ & 2.48\% & 17.3\% \\
Model 2 & Linear & 73,992 & $\delta=1,\ \gamma_D=0.6$ & 3.06\% & 21.9\% \\
\hline
Model 3 & Linear & 73,992 & $\delta=0.5,\ \gamma_D=0.85,\ \gamma_P=0.5$ & 2.56\% & 21.8\% \\
Model 3 & Linear & 73,992 & $\delta=0.25,\ \gamma_D=0.9,\ \gamma_P=0.5$ & 3.10\% & 31.1\% \\
\rowcolor{lightgreenrow}
Model 3 & Linear & 73,992 & $\delta=0.5,\ \gamma_D=0.9,\ \gamma_P=0.5$ & 2.26\%$\pm$0.27\% & 19.1\%$\pm$2.29\% \\
Model 3 & Linear & 73,992 & $\delta=0.75,\ \gamma_D=0.9,\ \gamma_P=0.5$ & 2.69\% & 22.9\% \\
Model 3 & Linear & 73,992 & $\delta=0.5,\ \gamma_D=0.95,\ \gamma_P=0.5$ & 3.09\% & 28.7\% \\
\hline
Model 4 & Nonlinear & 74,007 & $\delta=0.5,\ \gamma_{D_1}=\gamma_{D_2}=0.65,\ \gamma_{P_1}=\gamma_{P_2}=0.65$ & 3.81\% & 24.5\% \\
Model 4 & Nonlinear & 74,007 & $\delta=0.25,\ \gamma_{D_1}=\gamma_{D_2}=0.75,\ \gamma_{P_1}=\gamma_{P_2}=0.75$ & 2.89\% & 18.9\% \\
\rowcolor{lightgreenrow}
Model 4 & Nonlinear & 74,007 & $\delta=0.5,\ \gamma_{D_1}=\gamma_{D_2}=0.75,\ \gamma_{P_1}=\gamma_{P_2}=0.75$ & 2.35\%$\pm$0.05\% & 15.9\%$\pm$1.55\% \\
Model 4 & Nonlinear & 74,007 & $\delta=0.75,\ \gamma_{D_1}=\gamma_{D_2}=0.75,\ \gamma_{P_1}=\gamma_{P_2}=0.75$ & 2.92\% & 18.5\% \\
Model 4 & Nonlinear & 74,007 & $\delta=0.5,\ \gamma_{D_1}=\gamma_{D_2}=0.85,\ \gamma_{P_1}=\gamma_{P_2}=0.85$ & 3.98\% & 33.1\% \\
\hline
\end{tabular}

\caption{Architecture size, rheology, loss-weight sensitivity, and validation errors for the GNN models. Highlighted rows indicate the selected hyperparameter configurations used in the numerical experiments.}
\label{tab:model_hyperparams_sensitivity}
\end{table}

The GNN model employs the \emph{Gaussian Error Linear Unit} (GELU) as activation function and is trained using a processor composed of 30 message passing steps. 
Each message-passing layer includes 7 hidden layers, with skip connections introduced every 3 layers to enhance gradient propagation and model stability. 
Given these choices and considering the number of features, the number of trainable parameters is shown in Table \ref{tab:model_hyperparams_sensitivity}. The hyperparameters appearing in the loss function have been heuristically chosen to optimize the training process. A sensitivity analysis of such parameters is also reported in the same Table.

For Models~3--4, Table~\ref{tab:model_hyperparams_sensitivity} reports an exploratory hyperparameter sweep evaluated on the validation set. Each nonselected configuration was independently trained once, whereas the selected configuration for each model was retrained using five random seeds, for which the mean and standard deviation are reported. For both models, the five-seed mean errors of the selected configuration remain lower than the single-run errors of the neighbouring configurations for both pressure and velocity. Thus, the observed seed-to-seed variability does not materially change the hyperparameter selection within the tested range. Because the nonselected configurations were not evaluated across multiple seeds, however, the remaining differences should be interpreted as sensitivity trends rather than as statistically definitive rankings.

Each GNN model is trained using the Adam optimizer, with an initial learning rate of $10^{-3}$, decreased by a factor of 10 after 10 consecutive validation epochs without improvement. Training is stopped early if no improvement is observed for 25 epochs. We use a batch size of 1 (graph-level batching) and train for a maximum of 500 epochs.
GNN training and inference were performed on an NVIDIA GeForce RTX~4090. The linear and nonlinear full-order solves were run on an Intel Core i9-14900HX CPU; the workstation was equipped with 64~GB of RAM.
The timing entries in Tables~\ref{tab:mouse} and~\ref{tab:mousenonlin} compare the direct sparse full-order implementation with GNN inference on this workstation. They are implementation- and hardware-specific and are not an accelerator-neutral speed-up benchmark. A demonstration package is available at \texttt{\url{github.com/piermariovitullo/GNN-demo/}.}

\section{Numerical tests on GNN surrogate models}
\label{sec:results}
This section presents a detailed numerical investigation of the GNN-based models to evaluate their ability to accurately reproduce microvascular blood flow in a variety of four configurations and modeling complexities that are detailed in the following. The focus is on predicting the pressure, velocity and hematocrit distributions under linear and nonlinear rheological laws using the models of Section~\ref{sec:gnn_models}. 

\subsection{GNN tests on general-purpose capillary networks}

We now assess the generalizability of the four GNN models introduced in Section~\ref{sec:gnn_models} on a set of unseen vascular graphs, focusing on their ability to predict pressure, velocity and hematocrit fields beyond the training distribution. 

To quantify the accuracy of the model, we compute the mean relative $L^p$-error, $p=1,2$, with respect to the
high-fidelity (FOM) solutions for the predicted pressure, velocity and hematocrit fields. The $L^2$-error is consistent with the norm used in the loss function, and the $L^1$-error is a weaker norm for which lower error levels are expected.
For a generic output $\mathbf{U}^W=[\mathbf{U}^V,\mathbf{U}^E]^\top$ and for each testing graph $G^s_k$, we define
\begin{equation*}
    \varepsilon^*_k = 
    \frac{\| \mathbf{U}^*_k - \widehat{\mathbf{U}}^*_k \|_{p,*}}
         {\| \mathbf{U}^*_k \|_{p,*}} \times 100,
    \quad \text{and} \quad
    \overline{\varepsilon}^*= 
    \frac{1}{N_{\text{test}}} 
    \sum_{k=1}^{N_{\text{test}}} \varepsilon^*_k,
    \quad
    p=1,2, \ *=V,E, 
\end{equation*}
where $\widehat{\mathbf{U}}^W_k$ denotes the GNN prediction and $\overline{\varepsilon}^*$ the average error on all the graphs tested. 

\begin{figure}[t]
    \centering
    \includegraphics[height=0.2 \textheight]{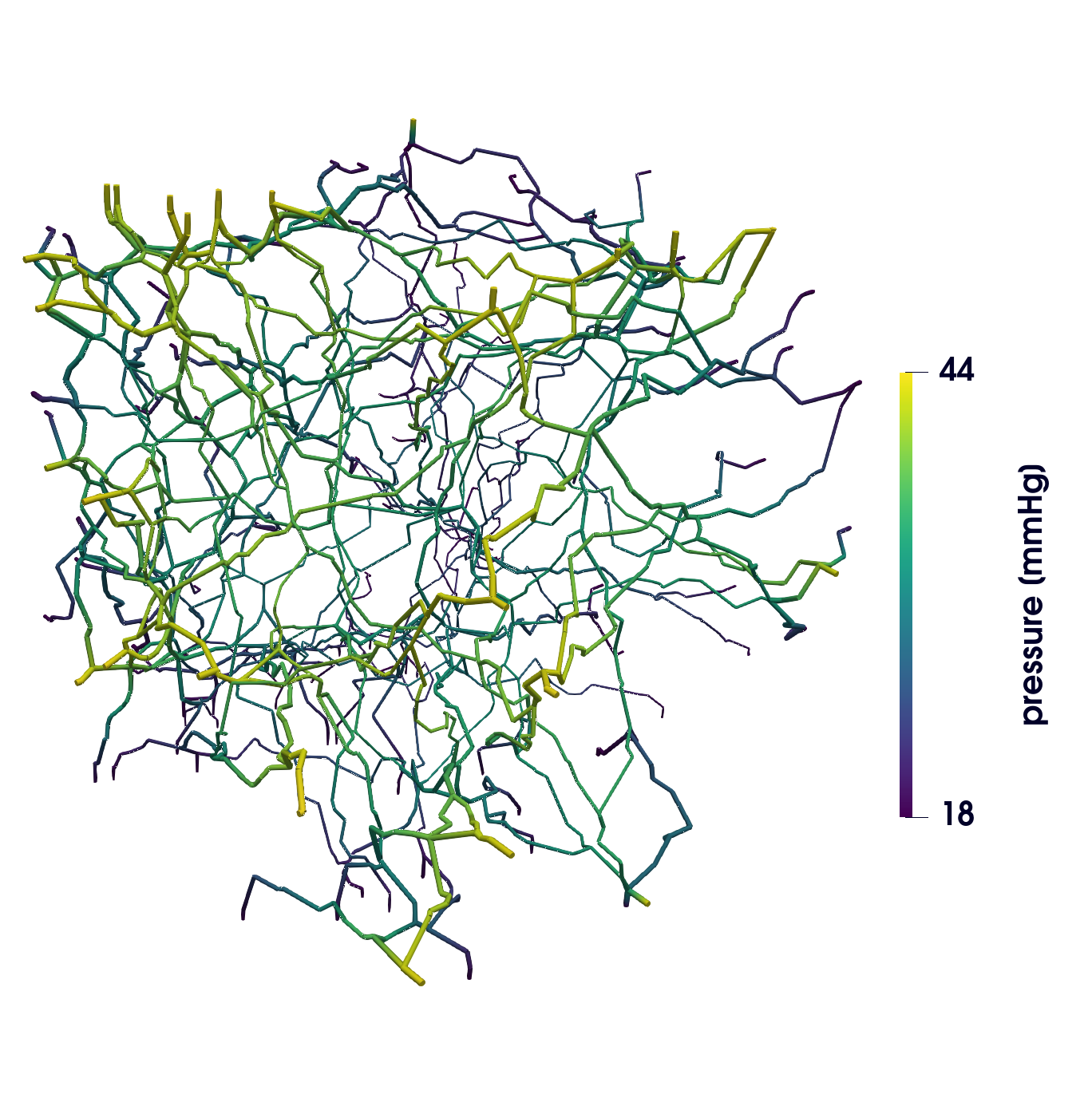}
    \includegraphics[height=0.2 \textheight]{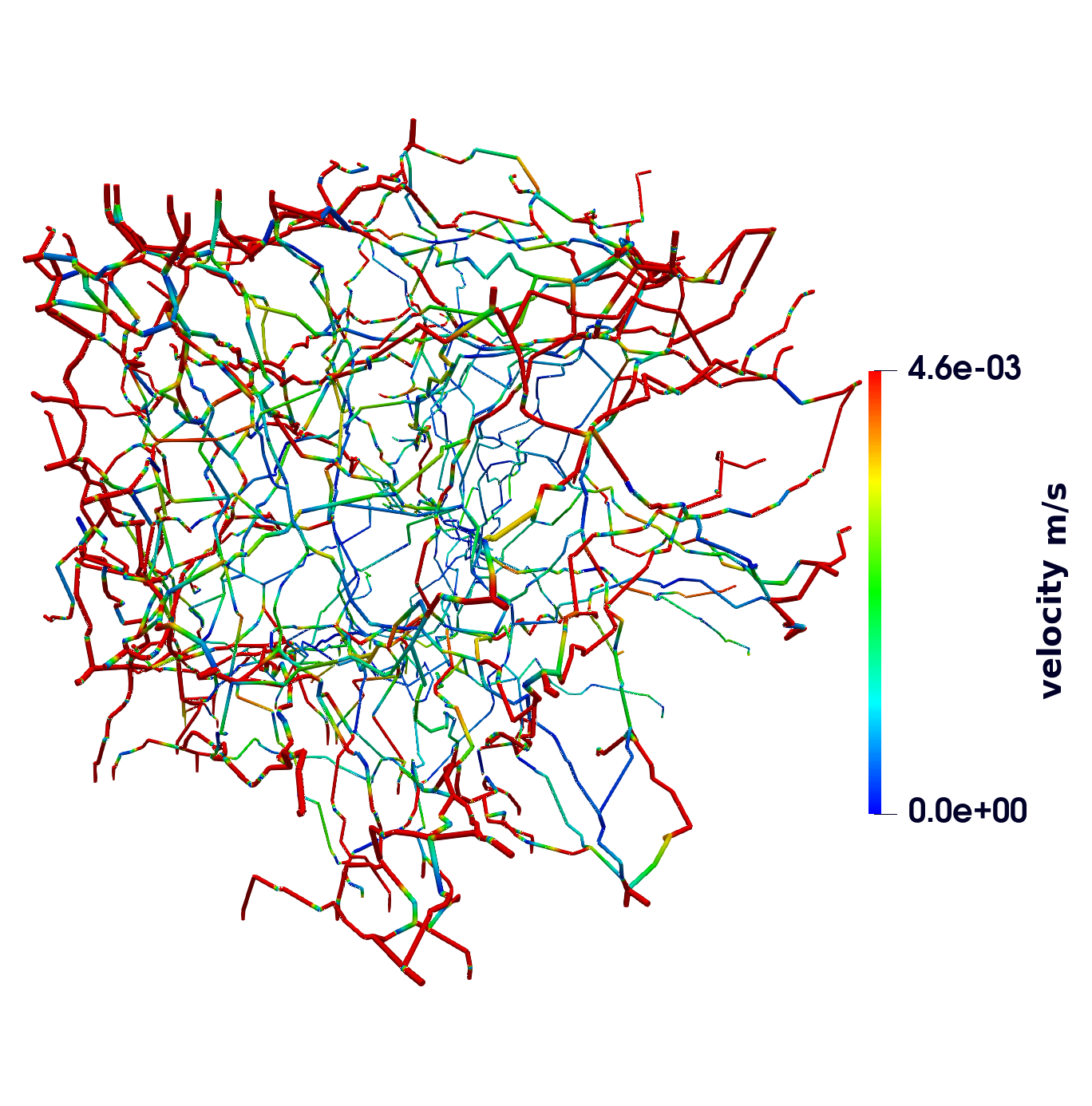}
    \includegraphics[height=0.2 \textheight]{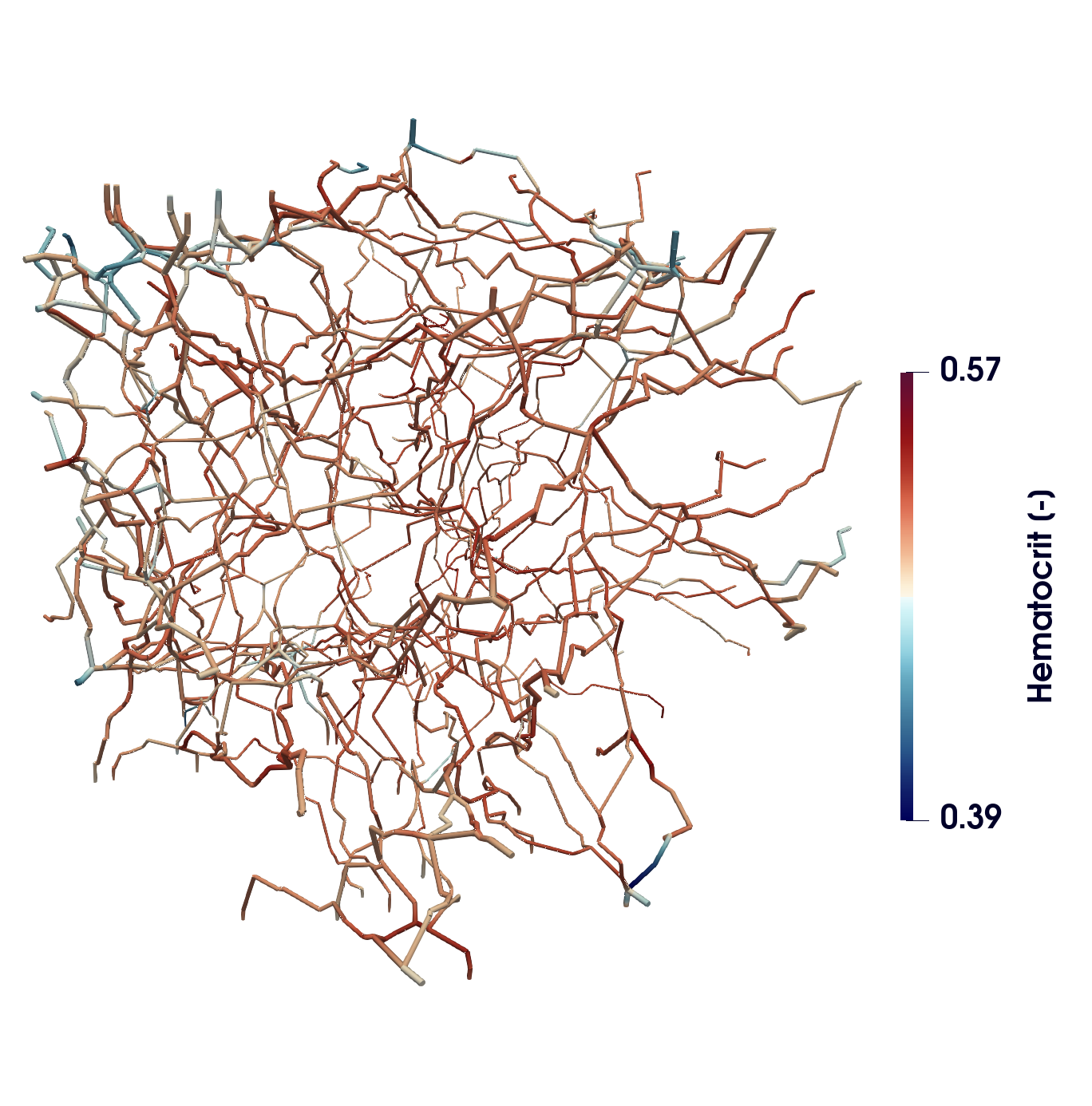} \\
    \vspace{-0.5cm}
    \includegraphics[height=0.2 \textheight]{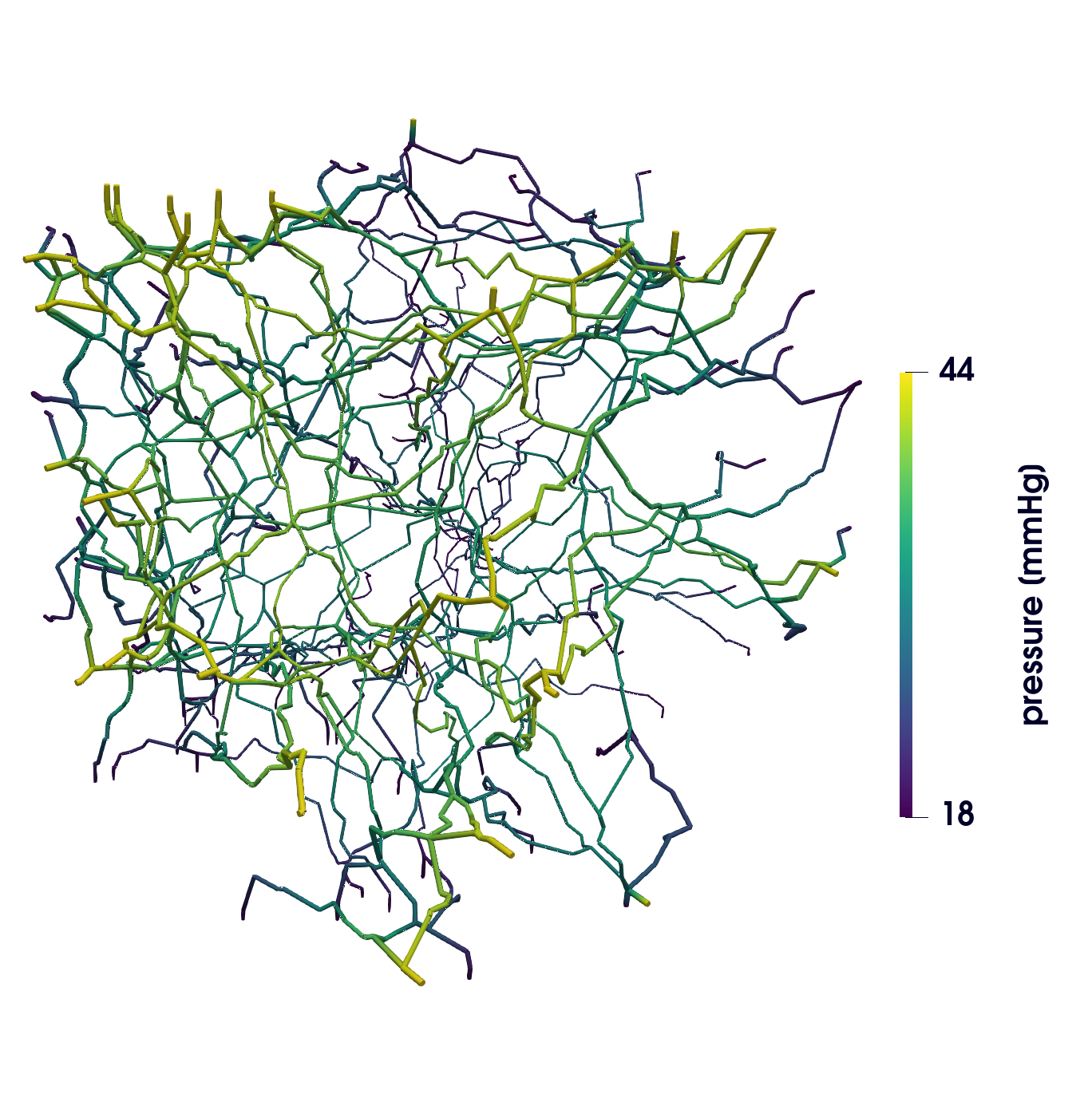} 
    \includegraphics[height=0.2 \textheight]{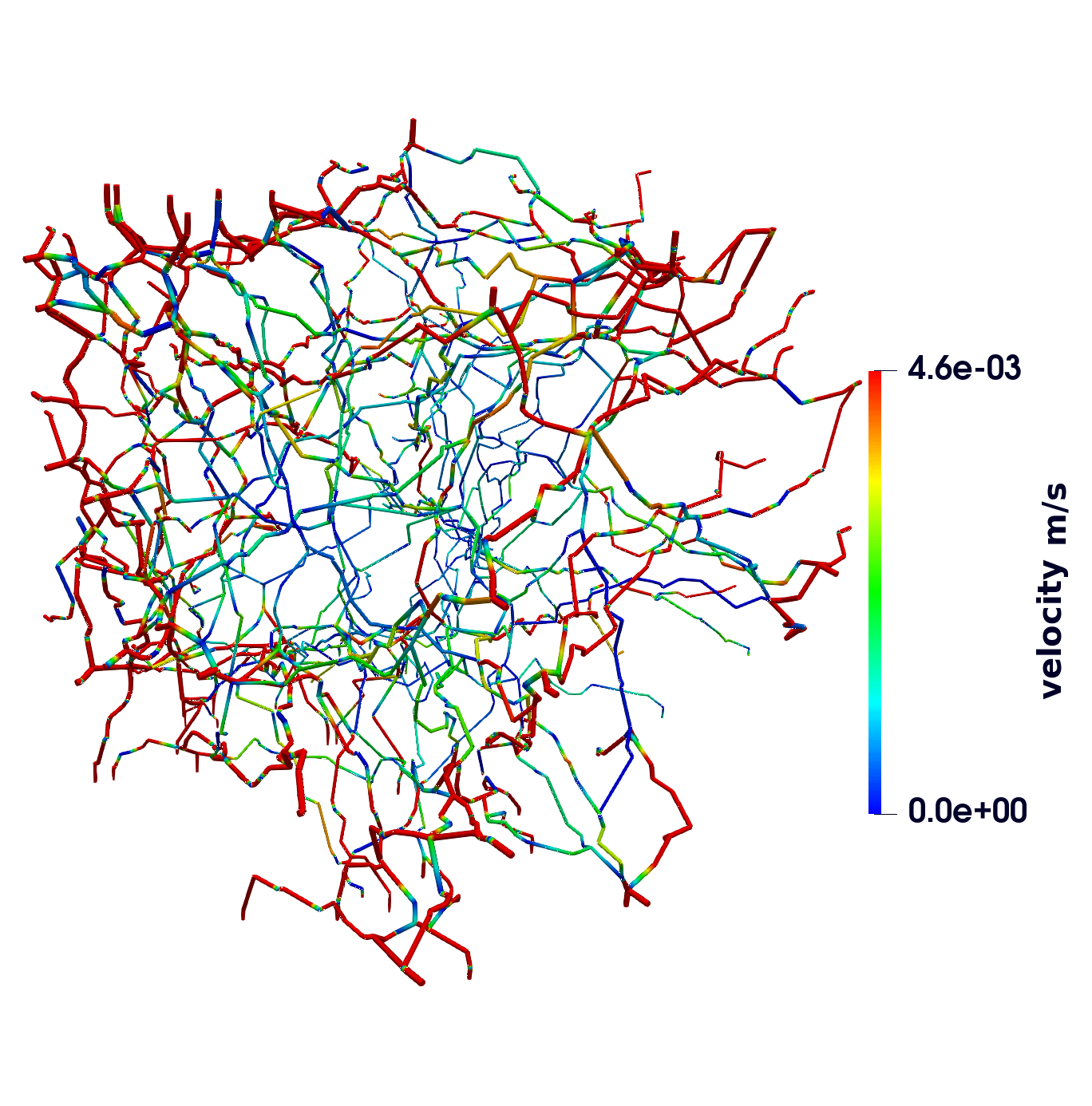} 
    \includegraphics[height=0.2 \textheight]{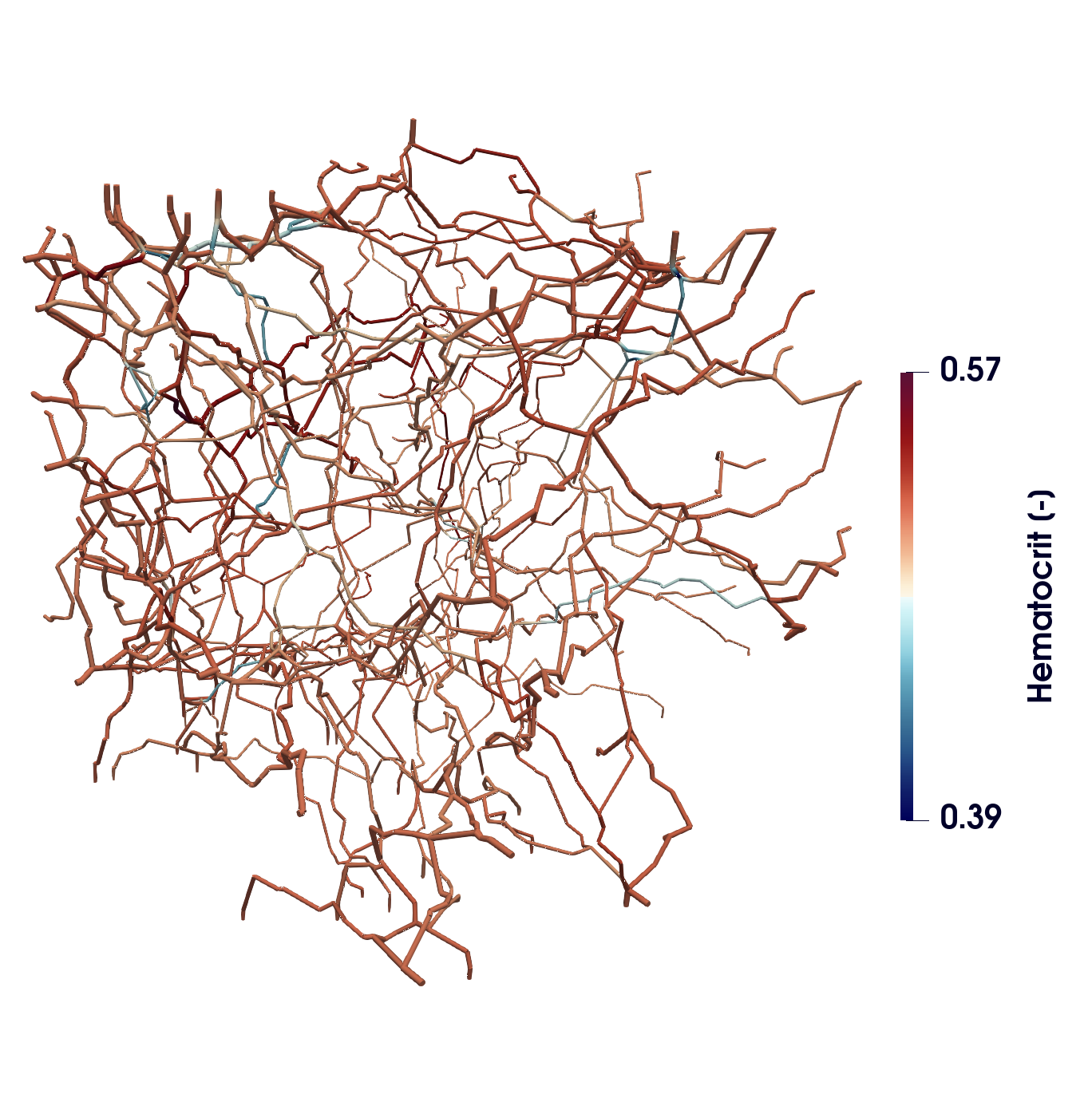}
    \caption[Visual comparison of pressure solutions for GNN model 4]{Visual comparison of pressure, velocity and hematocrit solutions (left, center, right panels, respectively) for GNN Model 4 (top), compared with the respective high-fidelity approximations (bottom). In this test case, the $L^{2}$ pressure, velocity and hematocrit GNN errors are respectively $2.55\%$, $16.02\%$ and $4.14\%$.}
    \label{fig:exp_visual_test}
\end{figure}

\begin{figure}[t]
\centering
\includegraphics[width=0.49\textwidth, trim=1mm 0.5mm 2mm 0.5mm, clip]{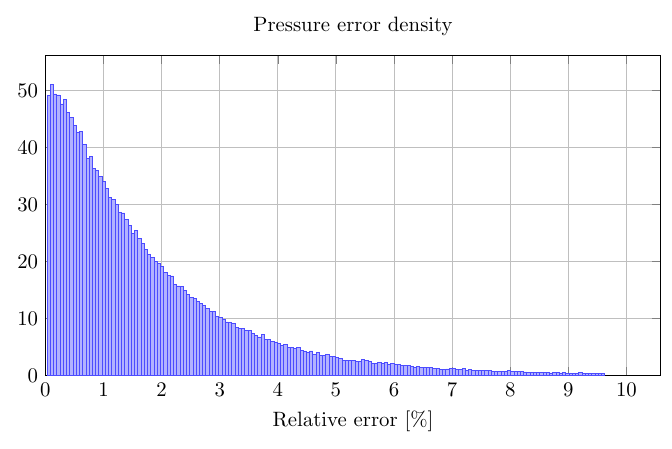}%
\hfill
\includegraphics[width=0.49\textwidth, trim=1mm 0.5mm 2mm 0.5mm, clip]{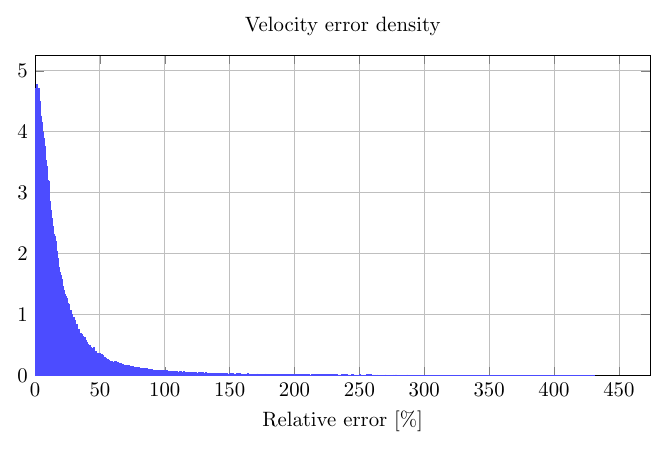}

\vspace{1mm}

\includegraphics[width=0.49\textwidth, trim=1mm 0.5mm 2mm 0.5mm, clip]{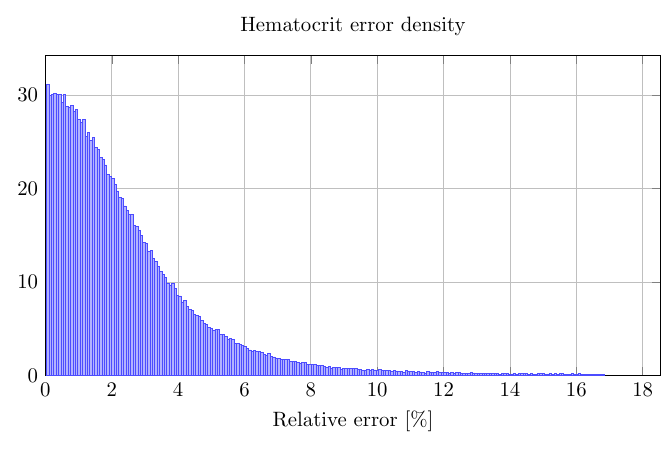}%
\hfill
\includegraphics[width=0.49\textwidth, trim=1mm 0.5mm 2mm 0.5mm, clip]{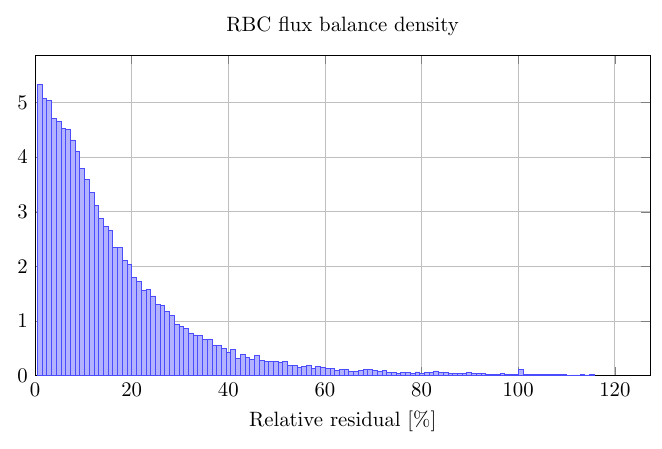}

\caption{Distributional error analysis (truncated at the 99th percentile for better visualization)} on the general-purpose vascular networks. Top row: empirical densities of the relative pressure and velocity errors. 
Bottom row: empirical density of the relative hematocrit error and of the relative RBC flux-balance residual at bifurcations.
\label{fig:error_postprocessing_pv}
\end{figure}

The corresponding values are summarized in Table~\ref{tab:err_combined} showing that the surrogate models comprehensively provide a good approximation.
For further insight, the visual comparison of the high-fidelity and surrogate solutions shown in Figure~\ref{fig:exp_visual_test} is not particularly informative. A quantitative analysis of the error distribution in the case of Model 4 is reported in Figure~\ref{fig:error_postprocessing_pv}, where we show the empirical distributions of pressure, velocity, hematocrit and RBC mass balance  errors.
The pressure error is strongly concentrated around small values and remains relatively well controlled across the network. By contrast, the velocity error displays a broader distribution, with a heavier tail associated with a smaller subset of vessel segments. 
The systematically larger velocity error reflects the different mathematical roles of pressure and velocity in the graph hemodynamic model. Pressure is a nodal potential governed by an elliptic, graph-Laplacian-like balance law; therefore, local prediction errors tend to be smoothed by the global network coupling and by the mass-conservation constraint. By contrast, velocity is an edge quantity derived from local pressure drops, resistances, and vessel diameters. Since resistance scales as \(D^{-4}\) and velocity scales as \(Q/D^2\), small errors in pressure differences, flow rates, or local geometric features may be amplified at the edge level, especially in small-diameter vessels. 

This interpretation is further supported by Supplementary Figure~\ref{fig:supp_velocity_flow_stratification}, where the per-edge relative velocity error is stratified by quantiles of the reference flow-rate magnitude $\lvert Q^{\mathrm{FOM}}\rvert$. Both the median error and its dispersion are largest in the lowest-flow quantiles and decrease as the flow rate increases, confirming that the heavy tail is primarily associated with low-flow vessel segments. The dependence on vessel diameter and topological distance from the nearest inlet is reported in Supplementary Figure~\ref{fig:supp_pressure_stratification}.
These error stratification data show that pressure error slightly increases in the small-diameter regime, where the solution is more sensitive to local geometric variations and to the stronger diameter dependence of the effective resistance.
In addition, the pressure error tends to increase with the topological distance from the nearest inlet.  Although the message-passing architecture exchanges information symmetrically across the graph, the target hemodynamic map remains inherently directional, being primarily organized by flow orientation from inlet to outlet. Note that this effect only appears in the general-purpose synthetic capillary networks, characterized by a uniform architecture that does not allow us to discriminate between arterioles, capillaries and venules.
In Figure~\ref{fig:error_postprocessing_pv}, the relative RBC mass-conservation residual evaluated at each bifurcation is also reported,
$$
r_j^{\mathrm{RBC}}
=
\frac{
\left|
\sum_{i\in E(v_j)} C_{ij} H_i Q_i
\right|
}{
\sum_{i\in E_{\mathrm{in}}(v_j)} H_i \lvert Q_i\rvert
},
$$
where the absolute RBC flux imbalance is normalized by the total incoming RBC flux. The figure reports \(100\,r_j^{\mathrm{RBC}}\) as a percentage. The resulting distribution is strongly right-skewed: most bifurcations are concentrated at relatively low residual values, while the density decreases steadily and develops a sparse tail extending beyond \(100\%\), revealing a limited number of bifurcations with substantial relative imbalance. 
The absolute RBC volumetric-flow imbalance is also concentrated at small values, with a median of $4.82\times10^{-14}\,\mathrm{m}^3\cdot\mathrm{s}^{-1}$ and a 95th percentile of $3.31\times10^{-13}\,\mathrm{m}^3\cdot\mathrm{s}^{-1}$, indicating that large absolute conservation violations are uncommon.

\begin{table}[t]
\centering
\setlength{\tabcolsep}{2pt}
\renewcommand{\arraystretch}{0.95}
\begin{tabular}{@{} |l|l|c|c|c| @{}}
\hline
\textbf{Model} & \textbf{Quantity} 
& \textbf{Training (L$^1$/L$^2$)} 
& \textbf{Validation (L$^1$/L$^2$)} 
& \textbf{Test (L$^1$/L$^2$)} \\
\hline
\multirow{3}{*}{\textbf{Model 1}}
& Pressure   & 2.11 / 2.95 & 2.19 / 2.94 & 2.20 / 2.83 \\
& Velocity   & -- / --   & -- / --   & -- / --  \\
& Hematocrit & --  / --     & --  / --     & --  / --     \\
\hline
\multirow{3}{*}{\textbf{Model 2}}
& Pressure   & 1.93 / 2.46 & 1.95 / 2.48 & 1.94 / 2.46 \\
& Velocity   & 12.5 / 17.0 & 12.7 / 17.3 & 12.6 / 17.7 \\
& Hematocrit & --  / --     & --  / --     & --  / --     \\
\hline
\multirow{3}{*}{\textbf{Model 3}}
& Pressure   & 1.66 / 2.21 & 1.72 / 2.29 & 1.65 / 2.20 \\
& Velocity   & 13.8 / 18.0 & 14.1 / 18.6 & 13.9 / 18.5 \\
& Hematocrit & --  / --     & --  / --     & --  / --     \\
\hline
\multirow{3}{*}{\textbf{Model 4}}
& Pressure   & 1.65 / 2.26 & 1.68 / 2.31 & 1.64 / 2.26 \\
& Velocity   & 13.5 / 13.3 & 13.7 / 13.9 & 13.7 / 14.1 \\
& Hematocrit & 2.59 / 3.75 & 2.60 / 3.76 & 2.59 / 3.75 \\
\hline
\end{tabular}
\caption{$L^1$ and $L^2$ relative errors (\%) with respect to the FOM solution for pressure, velocity, and hematocrit on the training, validation, and test sets. The reported values correspond to one representative random seed (the mean and standard deviation over multiple seeds are reported only in Table~\ref{tab:model_hyperparams_sensitivity}).}
\label{tab:err_combined}
\end{table}

\begin{figure}[t]
\centering
\includegraphics[width=0.49\textwidth, trim=1mm 0.5mm 2mm 0.5mm, clip]{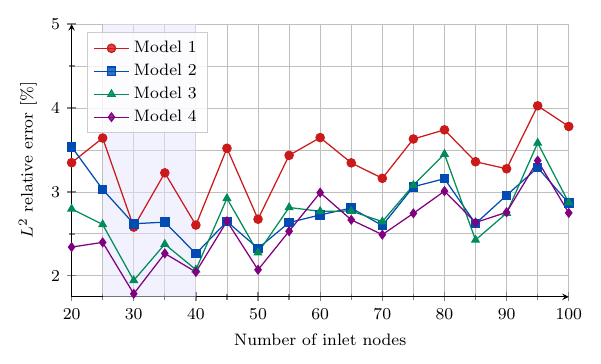}%
\hfill
\includegraphics[width=0.49\textwidth, trim=1mm 0.5mm 2mm 0.5mm, clip]{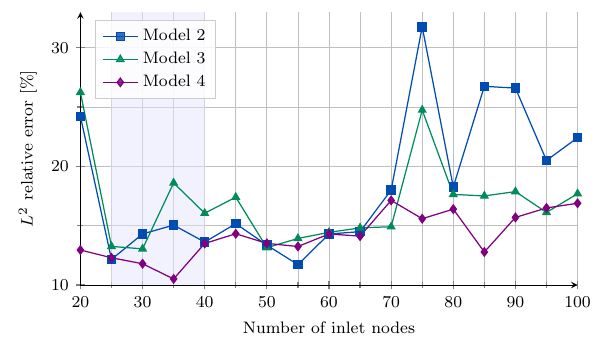}

\vspace{1mm}

\includegraphics[width=0.49\textwidth, trim=1mm 0.5mm 2mm 0.5mm, clip]{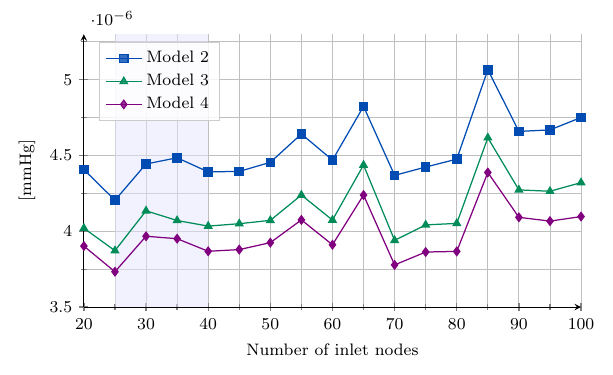}%
\hfill
\includegraphics[width=0.49\textwidth, trim=1mm 0.5mm 2mm 0.5mm, clip]{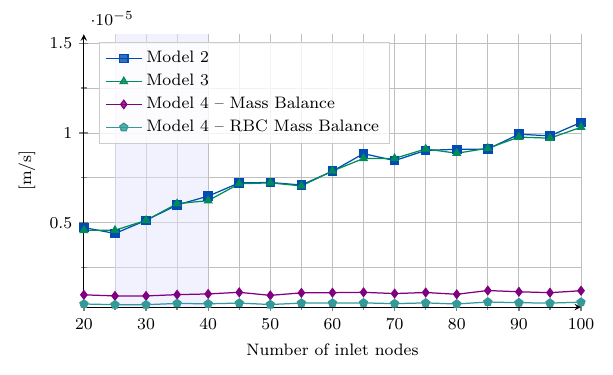}

\caption{Pressure (top-left) and velocity (top-right) $L^{2}$ relative error comparison between the different models. Physics-based residuals of the constitutive law (bottom-left) and the mass balance law (bottom-right). The gray band highlights the interval of number of inlet nodes where the training data are located. These data show the ability of the model to generalize well beyond the vascular complexity of the training dataset.}
\label{fig:error1}
\end{figure}

Figure~\ref{fig:error1} reports the $L^2$ mean relative pressure and velocity errors across the testing set for the four models (\(1\)--\(4\)). The evaluation is carried out in an \emph{out-of-sample} regime: the models were trained on vascular networks whose complexity, measured by the number of inlet nodes located on the upper surface of the domain, ranged between $25$ and $40$ (as indicated by the gray band in Figure~\ref{fig:error1}). 
The tests were then carried out on $17$ additional vascular graphs, characterized by a number of inlets varying from $20$ up to $100$, with each configuration containing twice as many outlet nodes as inlets.
The results highlight three main aspects. All GNN models exhibit stable accuracy far beyond the training window, effectively generalizing to vascular networks whose geometric complexity spans from $20$ to $100$ inlet nodes.  This holds consistently for pressure and velocity fields, confirming the robustness of the learned parameter-to-solution map $\mathcal{S}:(G,W)\mapsto \mathbf{U}^W$ to large variations in graph topology.
However, we note that, despite the increasing geometric complexity, the global hemodynamic regime remains coherent across all test graphs. This is confirmed by the near-constant depth of the networks—measured as the average number of nodes from inlets to outlets—whose weak dependence on the inlet count suggests that the flow regime is physically comparable across test configurations (not shown).
Finally, the physics-informed models (Model~3 and~4) maintain low residuals for constitutive and mass-balance laws throughout the entire testing range.  
More precisely, these metrics are defined as $\mathcal{L}_{\text{constitutive}}$ and $\mathcal{L}_{\text{mass}}$ for Model~3 and $\mathcal{L}_{\text{constitutive}}$, $\mathcal{L}_{\text{mass}}^{(1)}$, and $\mathcal{L}_{\text{mass}}^{(2)}$ for Model~4, in Section \ref{sec:gnn_models}.
These residuals remain small in all configurations, as shown in Figure~\ref{fig:error1} (bottom panel).
Overall, these results demonstrate that the proposed GNN architectures not only reproduce local flow quantities with high accuracy but also preserve the global physical structure of the problem, providing reliable predictions across a wide spectrum of vascular geometries.

\subsection{GNN tests on anatomically accurate cerebral vascular networks}

To evaluate the transfer of the proposed GNN pressure predictions, we test its performance on a set of anatomically realistic vascular networks generated using the methodology introduced in Section~\ref{sec:cco}. 
In contrast to the synthetic graphs used for training, these cerebral networks are significantly more complex, both topologically and morphologically.
In particular, physiological plausibility, hierarchical branching, and image-informed spatial constraints are now enforced.

We consider four synthetic networks representing the cortical circulation of the mouse brain, generated using algorithms developed in~\cite{linninger2,linningerplos,linninger1}. We name these graphs with the following labels: \texttt{S1.101}, \texttt{S2.102}, \texttt{S3.101}, and \texttt{S4.102}, and we refer to Table \ref{tab:morphometry} for a detailed description of their characteristics. These networks faithfully reproduce the statistical morphometric properties observed in two-photon laser scanning microscopy (2PLSM) datasets of mouse somatosensory cortex~\cite{Blinder2013,linninger2,linninger1}. Below, we provide a biological overview of their structure and functional roles.

\begin{table}[t]
    \centering
    \footnotesize
    \begin{tabular}{|l|c|c|c|c|}
        \hline
        \textbf{Network ID} & \textbf{S1.101} & \textbf{S2.102} & \textbf{S3.101} & \textbf{S4.102}  \\
        \hline
        \textbf{Number of nodes} & 66,090 & 88,286 & 151,626 & 154,466 \\
        \textbf{Number of segments} & 70,800 & 87,545 & 150,812 & 153,720 \\
        \textbf{Average diameter [$\mu$m]} & 5.72 & 5.85 & 5.67 & 5.59 \\
        \textbf{Mean segment length [$\mu$m]} & 40.1 & 39.6 & 41.0 & 38.8 \\
        \textbf{Capillary density [segments/mm$^3$]} & 11,300 & 11,500 & 11,600 & 11,400 \\
        \textbf{Total network length [mm]} & 89.6 & 113.4 & 209.8 & 211.7 \\
        \textbf{Total surface area [mm$^2$]} & 5.41 & 7.12 & 13.21 & 13.48 \\
        \textbf{Total volume [mm$^3$]} & 0.215 & 0.287 & 0.534 & 0.549 \\
        \hline
    \end{tabular}
    \caption[Morphometric summary of anatomically accurate vascular networks]{Morphometric properties of the four anatomically accurate mouse cortex networks used for the GNN generalization tests. These values are derived from the topological and geometric analysis of the synthesized graphs based on the CCO and microvascular closure algorithms.}
    \label{tab:morphometry}
\end{table}

Each vascular network represents an anatomically realistic reconstruction of the murine cortical microcirculation, encompassing interconnected pial arteries, penetrating arterioles, capillaries, and venules generated through the image-based Cerebral Network Synthesis (iCNS) framework~\cite{linningerplos}. The pial arteries distributed across the cortical surface serve as primary inflow routes, giving rise to $12$–$14$ penetrating arterioles per~mm$^2$, consistent with morphometric measurements~\cite{Blinder2013}. These vessels descend nearly orthogonally into the cortex, feeding a dense capillary plexus generated by a microvascular closure algorithm that guarantees topological continuity between the arterioles and venules. The capillary segments display diameter-dependent tortuosity (ranging from $1.1$ to $1.8$ for $d<150\,\mu$m), reproducing the geometric variability and flow dispersion observed in vivo. Downstream, post-capillary venules, and ascending veins collect blood and connect to pial veins that drain into the cortical sinuses, with venular densities between $13$ and $39$ per~mm$^2$, matching empirical data. Morphometric analysis confirms that segment counts ($\sim11{,}500$~mm$^{-3}$), as well as global quantities such as total vascular length, surface area, and volume, fall within the ranges reported for two-photon laser scanning microscopy (2PLSM) datasets of mouse somatosensory cortex~\cite{linningerplos}. 

Each network includes one inlet and one outlet node, representing the arterial and venous ends, respectively. These are the only two Dirichlet nodes used to compute the full-order reference solution. To construct the boundary features supplied to the GNN consistently with the training setup, we start from these nodes and identify the vertices connected to them through vessels with diameters greater than $12\,\mu$m. The pressures at the selected vertices are extracted from the full-order solution and supplied to the GNN as node features.
For Networks \texttt{S1.101}, \texttt{S2.102}, \texttt{S3.101}, and \texttt{S4.102}, the procedure selects, respectively, 156/56, 61/156, 361/396, and 761/571 inlet-side/outlet-side pressure-feature nodes. These correspond to $0.32\%$, $0.25\%$, $0.50\%$, and $0.86\%$ of all nodes in the respective networks. The selected nodes are used only as GNN input features and are not additional Dirichlet nodes of the full-order problem, which is solved with one inlet and one outlet. Although sparse relative to the complete network, these features are concentrated in the large-vessel compartments: in \texttt{S1.101}, for example, the 156 inlet-side and 56 outlet-side nodes represent $15.5\%$ of the arterial/arteriolar nodes and $7.0\%$ of the venous/venular nodes, respectively. Thus, the cortical test is conditional on a sparse set of reference-pressure values within these compartments, rather than on pressures supplied throughout the arterial and venous trees.

These anatomical networks provide physiologically consistent digital analogs of cerebral microcirculation, suitable for evaluating transfer of the GNN pressure predictions to realistic geometries.

\subsubsection{Tests on linear blood flow models}

Without retraining or fine-tuning the model, we apply GNN Model 3 to predict pressure fields over these domains. The results are summarized in Table~\ref{tab:mouse}.

Despite the large complexity gap and the anatomical specificity of the test domains, the surrogate model achieves an $L^2$-relative pressure error of less than 10.5\% in all cases. The reported implementation produces pressure predictions in under 30 milliseconds, compared with wall-clock times of about $0.4$--$1.3$ seconds for the direct sparse full-order implementation on the stated workstation. Because these timings are not a hardware-matched benchmark, they establish low-latency inference for the tested implementation rather than an implementation-independent speed-up. The pressure results show transfer to much larger, biologically realistic graphs, conditional on the supplied boundary-pressure features.

\begin{table}[t]
    \footnotesize
    \centering
    \begin{tabular}{|c c c c c|}
    \hline
     \textbf{Network} & \textbf{1} & \textbf{2} & \textbf{3} & \textbf{4}  \\
    \hline \hline
    \textbf{Pressure error ($L^1/L^2$)} & $5.06 \,/\, 6.34$ & $8.31 \,/\, 10.49$ & $7.64 \,/\, 8.80$ & $6.07 \,/\, 7.75$ \\
    \textbf{Nodes} & 66090 & 88286 & 151626 & 154466 \\
    \textbf{Solver time} & $0.43 s$ & $0.56 s$ & $1.13 s$ & $1.32 s$ \\
    \textbf{GNN time} & $23.32 ms$ & $27.01 ms$ & $28.39 ms$ & $26.58 ms$ \\
    \hline
    \end{tabular}
    \caption[Pressure error estimates for different complex vascular geometries]{$L^1$ and $L^2$ relative pressure error estimates (\%) for different anatomical vascular geometries, tested on the linear blood flow model compared with GNN model 3. The full-order columns report direct sparse-solver wall-clock time, while the GNN columns report inference time; the comparison is implementation- and hardware-specific.}
    \label{tab:mouse}
\end{table}

\subsubsection{Tests on nonlinear blood flow models}

We now assess its ability to generalize to nonlinear microvascular blood flow governed by hematocrit-dependent rheology. The resulting nonlinear system must be addressed using iterative techniques, i.e. the fixed-point iteration strategy described in Section~\ref{sec:nonlinear_models}, in which the nonlinear terms are treated using a nested loop approach. Precisely, the hematocrit field is held fixed while solving the linearized pressure and flow problem, and then updated via flow-dependent phase separation laws. Convergence is monitored by tracking changes in hematocrit between successive iterations.

\begin{figure}[t]
    \centering
    \includegraphics[height=0.21\textheight]{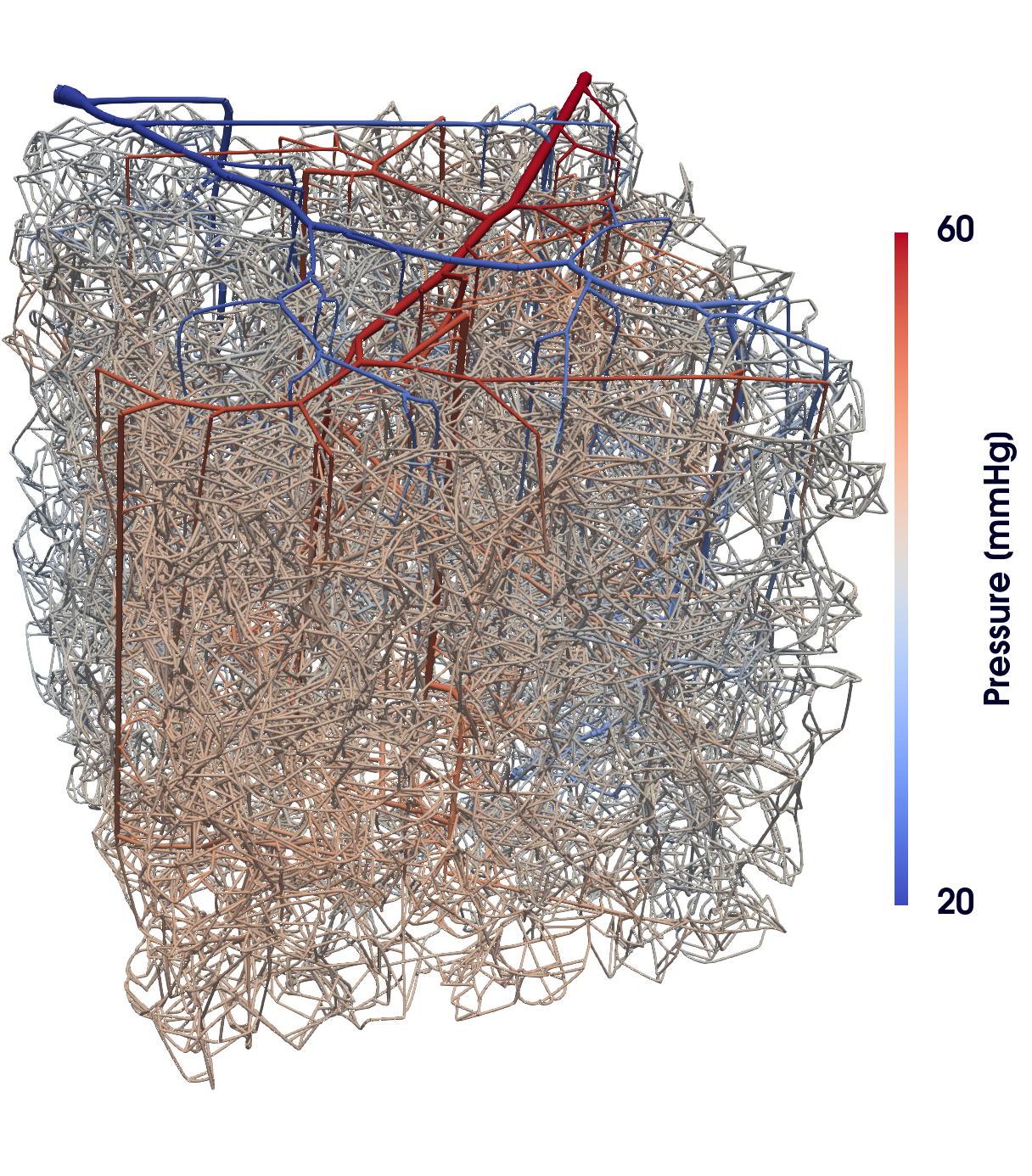}  
    \includegraphics[height=0.21\textheight]{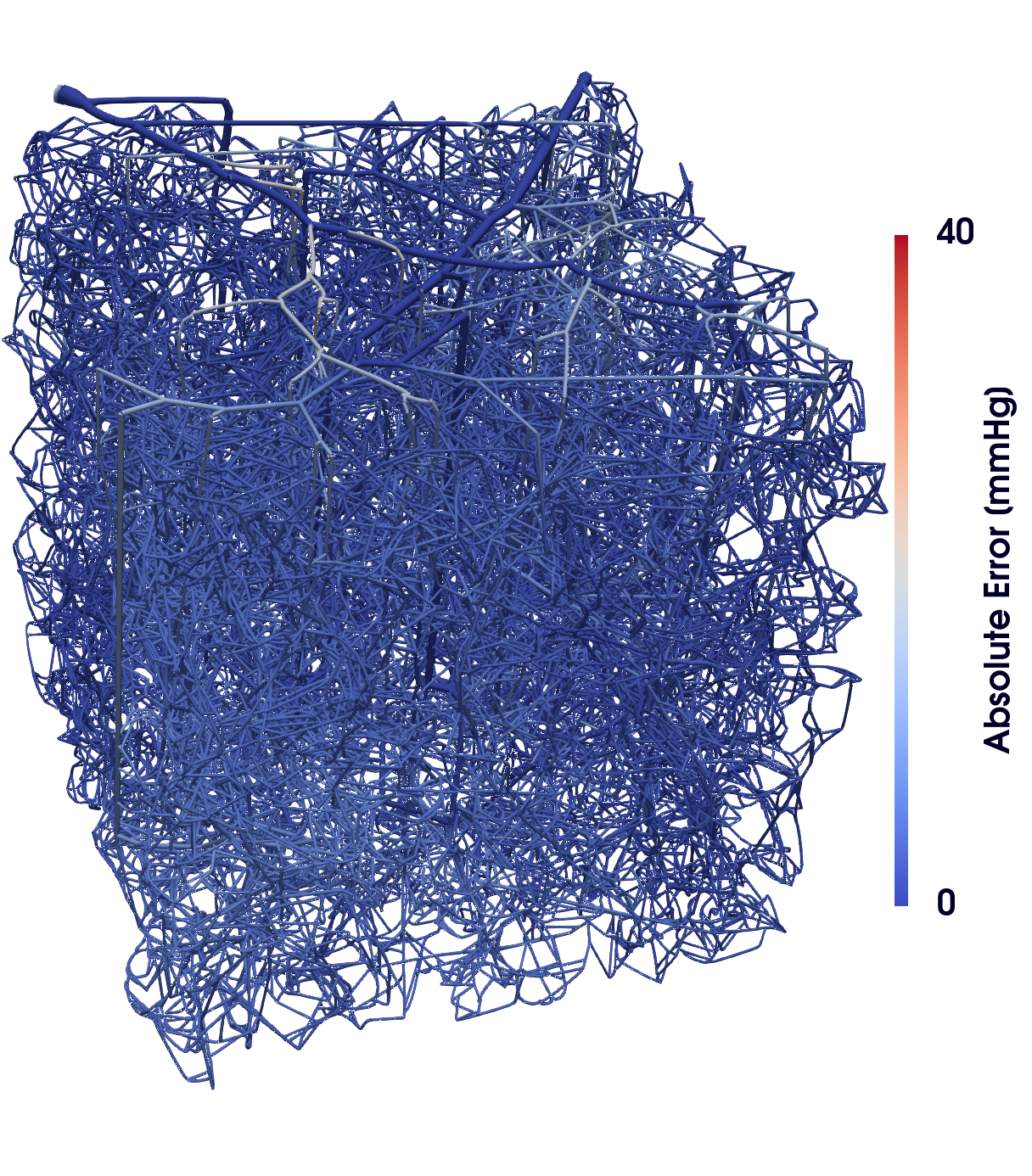}  
    \includegraphics[height=0.21\textheight]{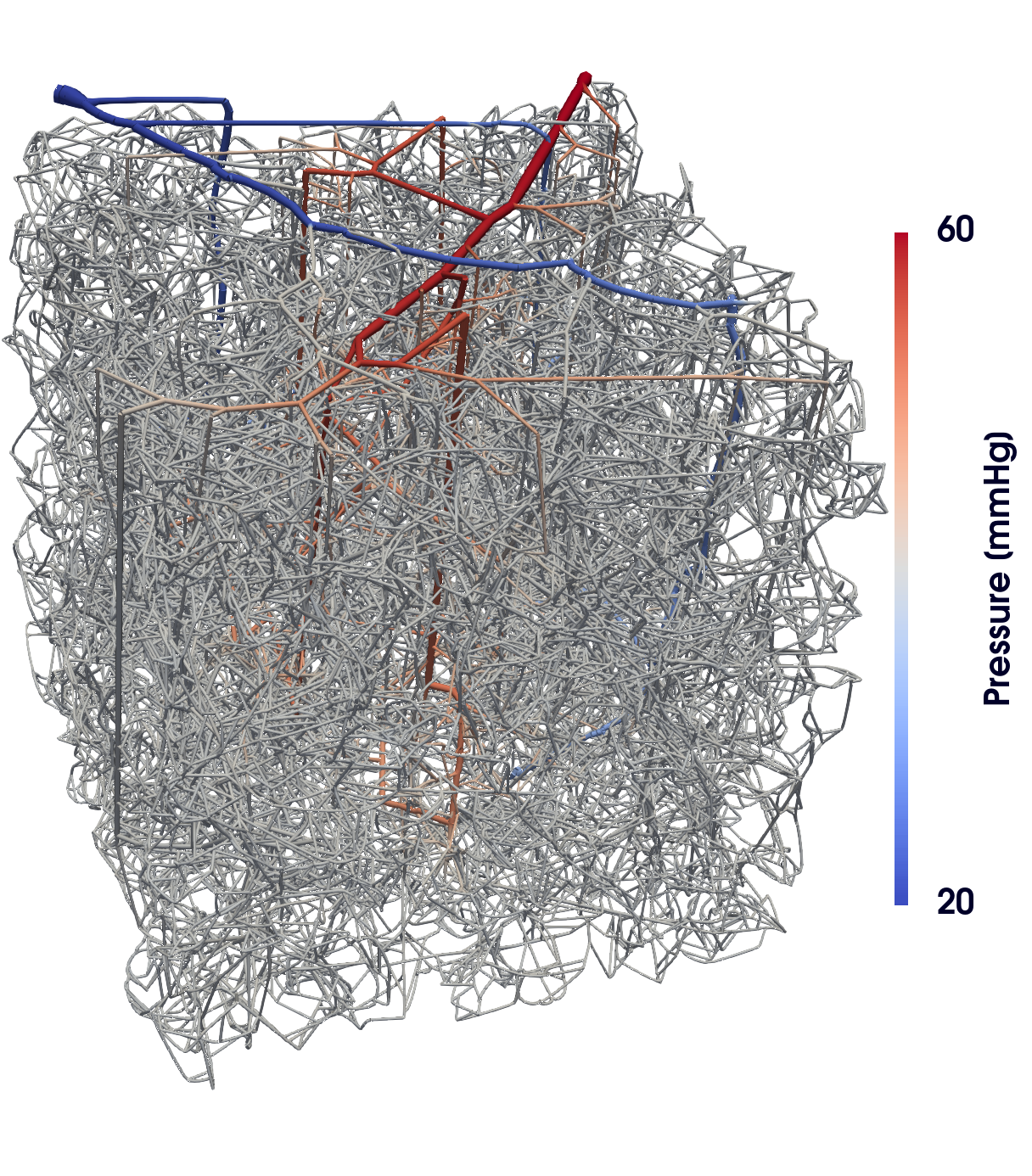}

    \includegraphics[height=0.35\textwidth, trim=1mm 0.5mm 2mm 0.5mm, clip]{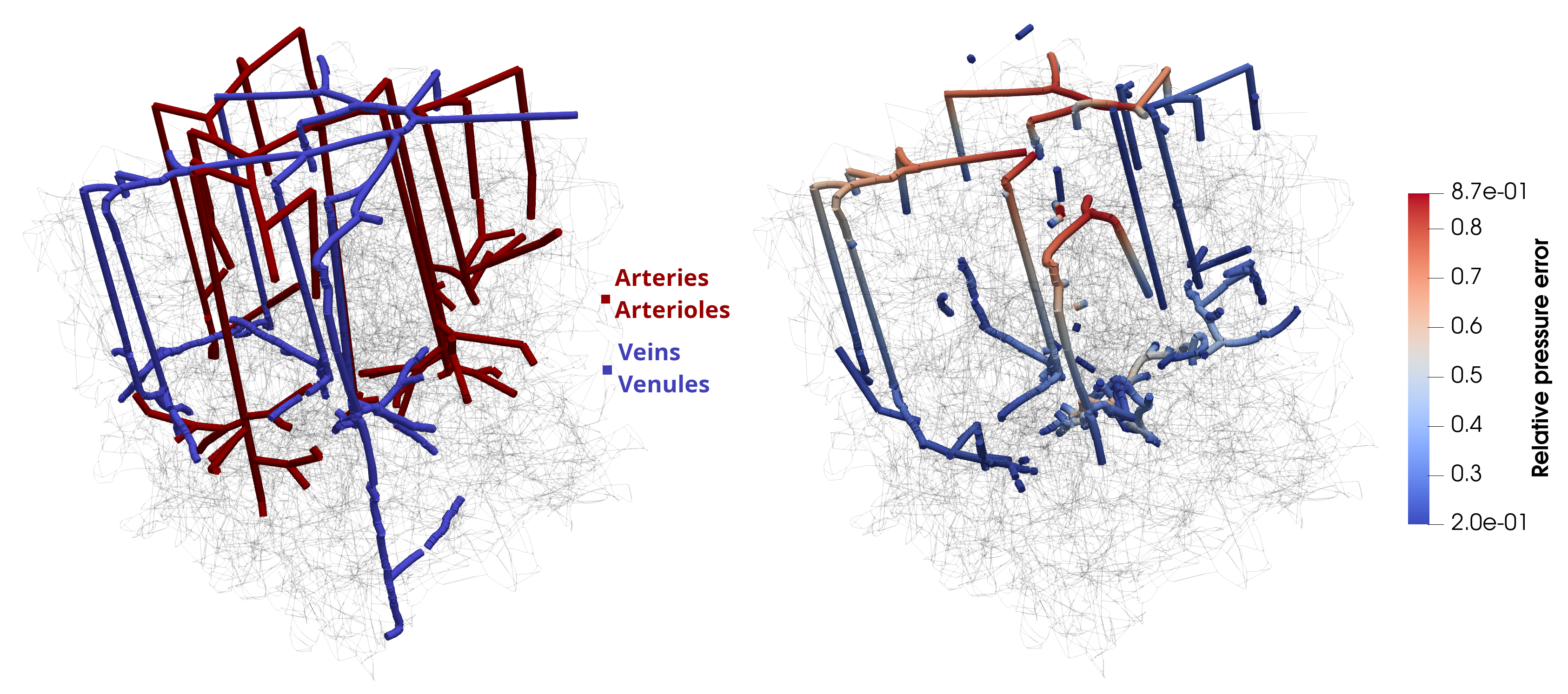}%
    \caption[Pressure prediction and spatial distribution of the error for Network \texttt{S1.101}]{Pressure prediction and spatial distribution of the error for Network \texttt{S1.101} in the nonlinear blood flow model. Top row: reference pressure field (left), absolute pressure error (center), and GNN-predicted pressure field (right). Bottom row: (left) visualization of the arterial/arteriolar (red)  and venous/venular compartments (blue), and spatial distribution of the largest relative pressure errors (right).}
    \label{fig:mouse3}
\end{figure}

In this test, we deploy GNN Model 4, trained exclusively on synthetic capillary networks with nonlinear rheology, on anatomically realistic cerebral networks. Although the architecture outputs nodal pressure, edge velocity, and hematocrit, the available cortical evaluation reports pressure only. We therefore restrict the quantitative generalization claim here to pressure and present those results in Figure~\ref{fig:mouse3} and Table~\ref{tab:mousenonlin}.

The surrogate model tested on nonlinear blood flow models achieves a $L^2$-relative pressure error of less than 9.5\% in all cases. GNN inference remains below 30 milliseconds in the reported implementation, whereas the fixed-point full-order runs require approximately $4$--$23$ seconds on the stated workstation. This is a hardware- and implementation-specific wall-clock comparison.

\begin{table}[t]
    \footnotesize
    \centering
    \begin{tabular}{|c c c c c|}
    \hline
     \textbf{Network} & \textbf{1} & \textbf{2} & \textbf{3} & \textbf{4}  \\
    \hline \hline
    \textbf{Pressure error ($L^1/L^2$)} & $6.08 \,/\, 8.06$ & $6.56 \,/\, 9.36$ & $6.09 \,/\, 9.14$ & $5.75 \,/\, 7.65$ \\
    \textbf{Nodes} & 66090 & 88286 & 151626 & 154466 \\
    \textbf{Solver time} & $4.29 s$ & $7.82 s$ & $13.52 s$ & $22.53 s$ \\
    \textbf{GNN time} & $23.32 ms$ & $22.01 ms$ & $28.35 ms$ & $29.51 ms$ \\
    \hline
    \end{tabular}
    \caption[Pressure error estimates for different complex vascular geometries]{$L^1$ and $L^2$ relative pressure error estimates (\%) for different complex vascular geometries, tested on the nonlinear blood flow model compared with GNN model 4. The full-order columns report fixed-point solver wall-clock time, while the GNN columns report inference time; the comparison is implementation- and hardware-specific.}
    \label{tab:mousenonlin}
\end{table}

\begin{figure}[t]
\centering
\includegraphics[width=0.47\textwidth, trim=1mm 0.5mm 2mm 0.5mm, clip]{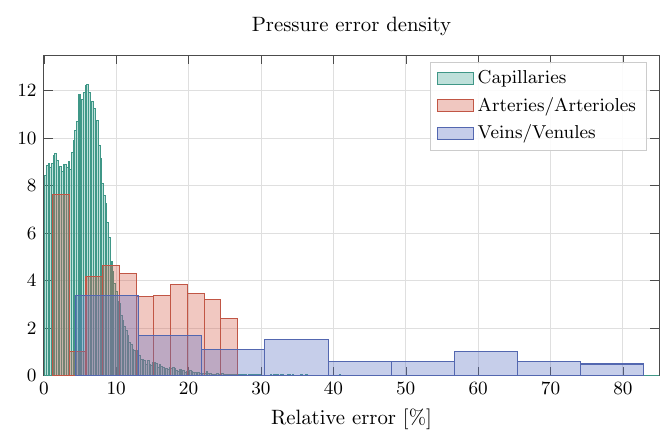}
\hfill
\includegraphics[width=0.49\textwidth, trim=1mm 0.5mm 2mm 0.5mm, clip]{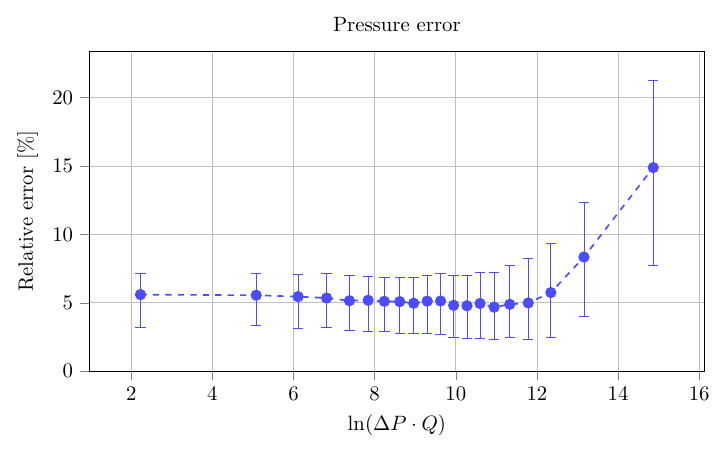}



\caption{Spatial and structural error analysis for pressure on cortical Network \texttt{S1.101}; the panels are not pooled across networks}.
\label{fig:mouse_pressure_postprocessing}
\end{figure}

Figure~\ref{fig:mouse3} (bottom row) suggests that the pressure error is not uniformly distributed across the anatomically realistic mouse networks, as the arteriolar and venular segments of the network systematically show errors larger than the average. A more detailed analysis is reported in Figure~\ref{fig:mouse_pressure_postprocessing} where vascular labels (Capillaries; Arteries/Arterioles; Veins/Venules) are used in the error distribution plot.
The left panel shows that errors are smallest in capillaries and larger in arterial and venous compartments, with the proximal arterial and distal venular sides being the most challenging regions.
For Network \texttt{S1.101}, without pooling data across different networks, the median and 95th-percentile relative pressure errors are, respectively, \(5.31\%\) and \(11.74\%\) in capillaries, \(11.43\%\) and \(25.18\%\) in arteries/arterioles, and \(23.97\%\) and \(78.14\%\) in veins/venules. These compartment-specific statistics quantitatively confirm that the pressure error is lowest in the capillary compartment and largest in the venous/venular compartment.
Consistently, the right panel of Figure~\ref{fig:mouse_pressure_postprocessing} shows that the error is directly associated with $\Delta P \cdot Q$, which represents the rate of mechanical energy dissipation along each vessel segment. The plotted abscissa is
\(
\ln\!\left(
\frac{\Delta P\cdot Q}
{1\,\mathrm{Pa}\,\mu\mathrm{m}^{3}\,\mathrm{s}^{-1}}
\right),
\)
that is, the natural logarithm of the numerical value of $\Delta P\cdot Q$ when $\Delta P$ is expressed in $\mathrm{Pa}$ and $Q$ in $\mu\mathrm{m}^{3}\,\mathrm{s}^{-1}$. The error increases sharply above a value of approximately \(12\) on this logarithmic scale. This indicates that the surrogate is less accurate in segments undergoing larger energetic loads, typically located near the inflow or outflow side of the network. By comparing these results with those obtained for generic capillary networks, we infer that, in anatomically realistic vascular networks, the error is mainly linked to the vessel hierarchy within the vascular tree. This occurs because the GNN models are attempting to predict microvascular flow under conditions where this hierarchical information is important, even though they were not trained with it.

\subsection{Discussion}

The results presented in this section demonstrate the ability of the proposed GNN framework to accurately approximate hemodynamic quantities on complex vascular networks while maintaining strong generalization properties across unseen topologies. 
Unlike traditional reduced-order or regression-based surrogates, which often rely on fixed connectivity or explicit graph parametrizations, the adopted message-passing formulation enables the model to operate directly on variable graph structures, making it inherently adaptable to different vascular architectures. Incorporating physics-informed loss terms further enhances the robustness of the surrogate model. 
By enforcing local mass conservation and constitutive relations linking pressure, velocity, and hematocrit, the proposed GNNs achieve not only low data-driven error but also strong physical consistency across a broad range of vascular configurations. 
The physics-informed terms regularize the learned solution toward the governing constraints, but the present comparison shows some mixed predictive effects. In fact, relative to Model~2, Model~3 improves pressure error slightly and worsens velocity error. The cortical experiments support pressure transfer to graphs with more than 150,000 segments; however, results about cortical velocity, hematocrit, and RBC-flux have not yet been evaluated on the cortical networks. Once trained, the GNN provides low-latency repeated evaluations of the parameter-to-solution map. This complements existing full-order whole-cortex solvers and may be useful in many-query workflows, although it is not an enabling factor to make a single steady whole-cortex solve feasible.


\section{Conclusions}\label{sec:conclusions}
This work introduces a physics-informed GNN framework for efficient simulation of microvascular blood flow. The surrogate model is trained on a large dataset of synthetically generated vascular graphs that capture generic capillary features. Its pressure predictions are subsequently evaluated on high-fidelity, anatomically accurate cerebrovascular networks representative of the mouse cortex.

We propose a two-tiered strategy that balances scalability and physiological realism. General-purpose vascular graphs generated via Voronoi tessellations provide a computationally inexpensive means of producing large, diverse training data. Anatomically constrained networks offer a challenging testbed to validate the surrogate’s extrapolation capabilities. These networks include complex topological features, heterogeneous vessel properties, and biological realism that cannot be achieved with heuristic generators. The GNN surrogate, trained solely on the former class of networks, provides evidence that its pressure predictions transfer to the latter, conditional on the supplied full-order-derived pressure features.
Our findings emphasize three key strengths of the proposed approach:  \textit{(i) Topology-aware modeling:} the graph-based architecture enables flexible learning across diverse vascular structures, supporting application to different organs and pathological conditions; \textit{(ii) Physics-informed learning:} incorporating governing-equation residuals constrains mass conservation and rheological consistency, while its effect on predictive error is quantity-dependent in the present experiments; \textit{(iii) Scalable inference:} the trained GNN surrogate provides low-latency prediction on large microvascular graphs, offering potential for repeated-query integration into multiscale simulation and digital-twin workflows.
Despite these strengths, current limitations include the restriction to steady-state flow and the assumption of impermeable, non-leaky vessels. Velocity and hematocrit predictions are also more sensitive to domain shifts, particularly in the presence of strong nonlinearities or unusual branching patterns.

Future work will focus on extending the learning framework introduced here in the context of \textit{capillary} flow (where blood motion is well described by a stationary regime dominated by diffusive mechanisms) to the case of arteriolar flow, where the regime is \textit{pulsatile}. In this case, the governing equations, the physical principles, and the resulting mathematical models differ substantially from those of capillary flow \cite{alastruey2011pulse,muller2016high}. For example, arterioles are accurately described by one-dimensional hyperbolic systems for pressure, flow, and cross/sectional area. The ongoing extension of the present methodology to the pulsatile setting is the subject of a dedicated study, whose preliminary results developed in \cite{behrens2025} are promising. Importantly, it confirms that the transition from capillaries to arterioles introduces substantial conceptual and numerical differences. Primarily, in pulsatile flow, the hierarchical organization of the vascular network becomes a key input feature in the learning architecture.
However, despite these differences in the underlying physics, the formulation of the computational learning problem remains conceptually analogous: in both cases, the hemodynamic state is represented on a graph, and a GNN can be designed to approximate the associated parameter-to-solution map. The message-passing mechanism, the encoder–processor–decoder architecture, and the treatment of nodal and edge variables remain essentially unchanged, demonstrating the flexibility and robustness of the GNN framework across different flow regimes. 

These observations reinforce the relevance and generality of graph-based surrogate models for hemodynamics, and show that the framework developed in this paper forms a solid foundation for future extensions to more complex vascular regimes.

\section*{Declarations}
\subsection*{Authorship and Contributorship Statement} All authors have made substantial intellectual contributions to the study conception, execution, and design of the work. All authors have read and approved the final manuscript.

\subsection*{Conflicts of Interest} The authors declare there are no conflicts of interest.

\subsection*{Data \& Code Availability Statement}

A binary demonstration package is publicly available at \texttt{\url{github.com/piermariovitullo/GNN-demo/}}. It provides two pretrained GNNs: the linear pressure--velocity model (Model~3, \texttt{EX1}) and the nonlinear pressure--velocity--haematocrit model (Model~4, \texttt{EX2}). Users can run inference on selected instances from three dataset groups: (\textit{i}) five synthetic Voronoi networks of increasing topological complexity; (\textit{ii}) 1200 synthetic Voronoi networks used for training, validation, and testing; and (\textit{iii}) one sample of the cortical mouse-brain network \texttt{S1}. Within this demonstration package, the network instances are accessed through compiled binary modules; the raw network files, the Python source code of the GNN modules, and the training code are not included. Separately, the \texttt{S1} cortical network is publicly available in ASCII format at \texttt{\url{github.com/LPPDAdmin/Microcirculation/}}. This repository provides the Cartesian coordinates of the graph nodes, edge connectivity and anatomical labels, vessel diameters, and a MATLAB demonstration for constructing the corresponding graph object. Thus, the publicly available resources support execution of the released pretrained models and direct access to the \texttt{S1} network geometry, but do not permit independent retraining or complete reproduction of all the results reported in this article.


\section*{Acknowledgments}
 PV and PZ acknowledge support from the MUR PRIN 2022 project 2022WKWZA8 \textit{Immersed methods for multiscale and multiphysics problems} (IMMEDIATE), part of the Next Generation EU program Mission 4, Component 2, CUP D53D23006010006; the European Union’s Euratom research and training programme 2021–2027 under grant agreement No. 101166699 \textit{Radiation safety through biological extended models and digital twins} (TETRIS);
 PZ acknowledges support from the \textit{Mathematical Challenges in Brain Mechanics} project of the Center for Advanced Study (CAS), The Norwegian Academy of Science and Letters, Oslo, Norway;
 PB acknowledges support from \textit{Dipartimento di Eccellenza} 2023--2027, Department of Mathematics, Politecnico di Milano. TV and AL gratefully acknowledge financial support through NIH NIA 1R01AG079894. PB, PV, PZ are members of the Gruppo Nazionale per il Calcolo Scientifico (GNCS) of the Istituto Nazionale di Alta Matematica (INdAM). 


\bibliographystyle{plain} 
\bibliography{biblio}


\appendix

\section*{Further Reading}

\begin{description}
    \item An introduction to microcirculation \cite{Popel2005};
    \item Synthesis and challenges in the simulation of large scale vascular networks \cite{linninger1,linninger2};
    \item An introduction to learning algorithms and scientific machine learning \cite{doi:10.1142/S0218202525500125};
    \item An introduction to Graph Neural Networks \cite{Bronstein2017}; 
\end{description}

\newpage

\section*{Supplementary Material}

\subsection*{Pressure-error stratification on general-purpose networks}

The analysis below uses Model~4 on the 120-graph test set. The data are partitioned into 20 quantile bins of equal cardinality, each containing 5\% of the overall samples. Each point is a bin median and the bars span the interquartile range. For a pressure node, the diameter coordinate is the mean diameter of its incident edges. Topological distance is the shortest number of graph edges from the nearest inlet to the node at which pressure error is evaluated.

\begin{figure}[H]
\centering
\includegraphics[width=0.49\textwidth, trim=1mm 0.5mm 2mm 0.5mm, clip]{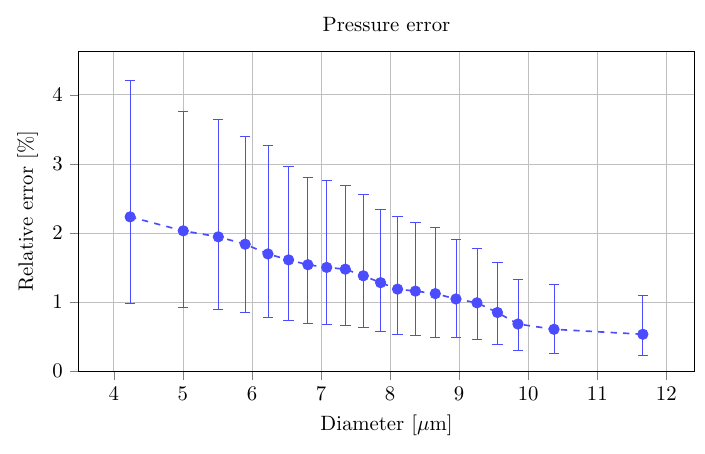}%
\hfill
\includegraphics[width=0.49\textwidth, trim=1mm 0.5mm 2mm 0.5mm, clip]{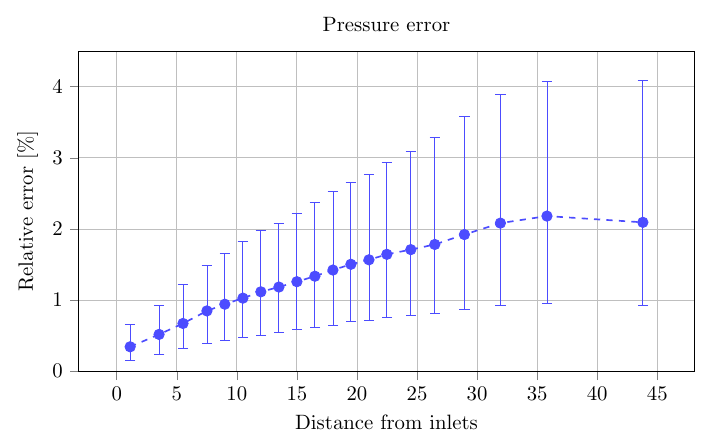}
\caption{Pressure relative error stratified by vessel diameter (left) and topological distance from the nearest inlet (right) on the general-purpose test networks. The medians increase toward small diameters and with distance from the nearest inlet; the broad interquartile ranges show substantial within-bin variability.}
\label{fig:supp_pressure_stratification}
\end{figure}

\subsection*{Velocity-error stratification by flow rate}

To distinguish the behavior of the per-edge relative error from that of the global \(L^2\) metric, we stratify the velocity error according to quantiles of the reference flow-rate magnitude \(\lvert Q_i^{\mathrm{FOM}}\rvert\). The analysis uses Model~4 on the 120-graph test set, and the velocity error is evaluated after recovering the physical axial velocity. The data are partitioned into 20 quantile bins of equal cardinality, each containing 5\% of the overall samples. Each point represents the median error within a flow-quantile bin, its horizontal coordinate is the mean reference flow magnitude in that bin, and the bars span the interquartile range.

\begin{figure}[H]
\centering
\includegraphics[width=0.65\textwidth,trim=1mm 0.5mm 2mm 0.5mm,clip]{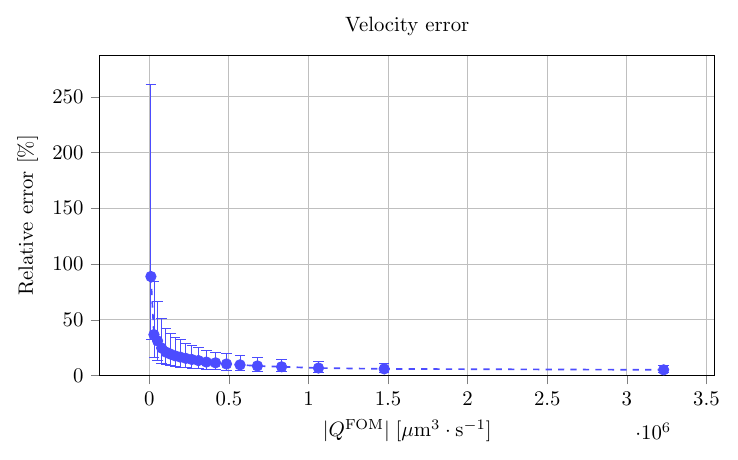}
\caption{Per-edge relative velocity error stratified by quantiles of the reference flow-rate magnitude $\lvert Q^{\mathrm{FOM}}\rvert$ on the general-purpose test networks. The largest relative errors and the broadest interquartile range occur in the lowest-flow quantiles, whereas both the median error and its variability decrease rapidly as the flow rate increases.}
\label{fig:supp_velocity_flow_stratification}
\end{figure}



\end{document}